%% file: HP-s0.tex
\input amstex
\documentstyle{amsppt}
\magnification=1200
\hsize=13.8cm
\vsize=19cm
\catcode`\@=11
\def\NoLogo{\let\logo@\empty}
\catcode`\@=\active
\NoLogo

\def\heat{\lf(\frac{\p}{\p t}-\Delta\ri)}
\def \b {\beta}
\def\i{\sqrt{-1}}
\def\Ric{\text{Ric}}
\def\lf{\left}
\def\ri{\right}
\def\bbar{\bar \beta}
\def\a{\alpha}
\def\g{\gamma}
\def\p{\partial}
\def\delbar{{\bar\delta}}
\def\ddbar{\partial\bar\partial}

\def\C{\Bbb C}

\def\tD{\tilde \Delta}

\def\ba{{\bar\alpha}}
\def\bb{{\bar\beta}}
\def\abb{{\alpha\bar\beta}}

\def\i{\sqrt {-1}}
\def\tD{\widetilde \Delta}
\def\tn{\widetilde \nabla}
\def \D {\Delta}
\def\aint{\frac{\ \ }{\ \ }{\hskip -0.4cm}\int}
\documentstyle{amsppt}
\magnification=1200
\hsize=13.8cm

\leftheadtext{Lei Ni  and Luen-fai Tam}
\rightheadtext{Plurisubharmonic functions, K\"ahler-Ricci flow}
\topmatter
\title{Plurisubharmonic functions and  the K\"ahler-Ricci flow}\endtitle

\author{Lei Ni\footnotemark and  Luen-Fai Tam\footnotemark}\endauthor
\footnotetext"$^{1}$"{Research partially supported by NSF grant
DMS-0196405 and DMS-0203023, USA.}
\footnotetext"$^2$"{Research partially supported by   Earmarked Grant of
Hong Kong \#CUHK4217/99P.}
\address
Department of Mathematics, University of California, La Jolla, CA 92093
\endaddress
\email{
lni\@math.ucsd.edu}
\endemail

\address
Department of Mathematics, The Chinese University of Hong Kong,
Shatin, Hong Kong, China
\endaddress
\email{lftam\@math.cuhk.edu.hk}
\endemail

\affil
{
Stanford University\\
The Chinese University of Hong Kong
}
\endaffil

\date  October, 2001
\enddate

\endtopmatter

\document

\subheading{\S0 Introduction}
Consider the K\"ahler-Ricci flow on a K\"ahler manifold $(M,g_\abb(x))$:
$$
\frac{\p}{\p t}g_\abb=-R_\abb,\ \ g_{\abb}(x,0)=g_{\abb}(x) \tag0.1
$$
In this work, $(M,g_\abb(x))$ will be assumed to be complete and noncompact
with bounded nonnegative holomorphic bisectional curvature. Solutions
 of (0.1)
on complete noncompact K\"ahler manifolds with bounded nonnegative
bisectional
 curvature was extensively studied in a series of paper of Shi \cite{Sh1-3}.

Important properties and applications have also been obtained, see
 \cite{Sh2-3, C-Z, C-T-Z}.  In [N-T], the authors studied the
 K\"ahler-Ricci
flow from another point of view. Namely, solutions of (0.1) are
 investigated by
 using the solution to
the Poincar\'e-Lelong equation obtained in \cite{M-S-Y, N-S-T1}.
More precise,
it was proved in \cite{N-S-T1}, under
some mild average assumptions on scalar curvature ${\Cal R}_0(x)$ of
the initial
 metric $g_\abb(x)$, one can solve
the Poincar\'e-Lelong equation:
$$
\i\ddbar u_0=\Ric_0,\tag0.2
$$
where $\Ric_0$ is the Ricci form of the initial metric. We should mention
that
 (0.2) was solved by Mok-Siu-Yau \cite{M-S-Y} and Mok \cite{M1} in case
when
$M$ has maximal volume growth and the scalar curvature  has  quadratical
pointwise decay. Using the solution of (0.2),   one can easily find a
function
 $u(x,t)$ so that $\i\ddbar u(x,t) =\Ric(x,t)$ where $\Ric $ is the Ricci
 form
of the metric $g(t)$. Moreover $u(x,t)$ satisfies  time-dependent
heat equation:
$$
\cases&\heat u=0\\
&u(x,0)=u_0(x).
\endcases\tag0.3
$$
Hence one  can study (0.1) by using (0.3). For example, we gave a
simple proof of the long time existence of (0.1) under some growth
conditions on
${\Cal R}_0$ in [N-T].  Note that in this case $u_0$ and
 $u(\cdot,t)$ are both
plurisubharmonic because $g(t)$ has nonnegative holomorphic bisectional
curvature \cite{Sh2-3}.

This motivates us to  study  (0.3) for general plurisubharmonic initial
 data
$u_0$. Under a rather mild assumption  on the growth rate of $u_0$, we
 can prove
 that (0.3) has a long time solution in the sense that if the
 K\"ahler-Ricci
 flow (0.1) has a solution up to time $T$,
then (0.3) also has a unique solution
up to time $T$.

The next  important question is whether or not under the flow (0.3) the
plurisubharmonicity will be preserved. In order to study this problem, we
have to study the complex Hessian $u_\abb$ of $u$. One can show that
 $u_\abb$
satisfies the complex Lichnerowicz-Laplacian heat equation (see (1.2)
for details). We shall prove that if the $\tD u_0$ is of at most
exponential
growth, then plurisubharmonicity will be preserved. Here $\tD$ is
the Laplacian
 of the initial metric. In fact, we shall prove the result for
 more general
Hermitian symmetric (1,1) tensors which satisfy the
Lichnerowicz-Laplacian heat
 equation. See Proposition 1.1. In case $u_0$ is the solution of (0.2),
 the assumption on the rate of growth of $\tD u_0$ is the same as the
assumption
 on the rate of growth of the
scalar curvature $\Cal R_0$, which is assumed to
 be bounded in \cite{Sh2-3}.

There are many important differential Harnack type
inequalities for Ricci flow
and curvature flows obtained by various people, see
\cite{L-Y, H4, Cw1-2, Co1-2, A} for examples. Works in this area can be
traced back  to the fundamental works of Li-Yau \cite{L-Y} and Hamilton
\cite{H4}.  For this reason, in this paper we shall call this kind of
inequalities to be Li-Yau-Hamilton type inequalities, or LYH  inequalities
for short.
 In \cite{C-H}, Chow and Hamilton  obtained a linear trace LYH
 inequality
for a symmetric two-tensor on a Riemannian manifold with a family of metric
$g(t)$ satisfying the Ricci flow equation so that the initial metric has
nonnegative curvature operator. The two-tensor is assumed to satisfy the
real Lichnerowicz-Laplacian heat equation.  In  this paper, using the
results
 of Cao \cite{Co1-2}, we shall prove a complex version of
 Chow-Hamilton's result.
 More precisely, suppose (0.1) has a solution on $M\times[0,T]$ so that the
 initial metric has nonnegative bounded holomorphic bisectional curvature.
 Let $h_\abb(x,t)$ be a solution of the complex Lichnerowicz-Laplacian heat
equation so that $h_\abb(x,0)\ge 0$ and $h_\abb(x,t)$ satisfies some growth
conditions. Then on $M\times (0,T]$, we have
$$
\split
Z& =\frac12[g^{\a\bbar}\nabla_{\bbar} div(h)_\a+g^{\g\delbar}
\nabla_\g div(h)_\delbar]\\
&\quad+g^{\a\bbar}g^{\g\delbar}[R_{\a\delbar}h_{\g\bbar}+
\nabla_\g h_{\a\delbar}V_{\bbar}+
\nabla_{\bbar}h_{\a\delbar}V_{\g}+h_{\a\delbar}V_{\bbar}V_\g ]+
\frac{H}{t}\ge 0\endsplit
$$
where   $div(h)_\a= g^{\g\delbar}\nabla_\g h_{\a\delbar}$,
 $div(h)_{\bbar}=g^{\g\delbar}\nabla_{\delbar}h_{\g \bbar}$,
$H$ is the trace of $h_{\a\bbar}$ with respect to $g_{\a\bbar}(x,t)$,
and $V$ is any vector field of type (1,0). See Theorem 1.1 for details.
 In case $h_\abb$ is the complex Hessian of a plurisubharmonic solution
 of (0.3) so that $u_0$ is not harmonic, then our result implies that
$w=u_t$ satisfies
$$
w_t-\frac{|\nabla w|^2}w+\frac wt\ge0
$$
which   extends Cao's trace LYH type inequality for the scalar curvature \cite{Co1-2}. Unlike \cite{C-H}, which mainly considers compact manifolds and not very specific on noncompact manifolds, we need the growth conditions on $h_\abb$ so that one can apply the maximum principle in \cite{N-T}.

As an application of the above results on the study of (0.3) and the linear trace LYH type inequality, we shall study Liouville properties for plurisubharmonic functions on $(M,g_\abb(x))$. Suppose (0.1) has long time solution. Then we have the following:

\proclaim{Theorem 3.1} Let $(M,g_\abb(x))$ be a complete noncompact manifold with bounded nonnegative holomorphic bisectional curvature so that the (0.1) has long time solution. Suppose $u_0$ is a plurisubharmonic function such that (i) $u$ is bounded; and (ii) $\tD u_0(x)\le \exp(a(1+r_0(x))$ for some constant $a>0$. Then $u_0$ must be constant.
\endproclaim

In $\C^m$, a plurisubharmonic function with sub-logarithmic growth must be
 constant. It is conjectured that this is still true for
 complete noncompact K\"ahler manifolds with nonnegative Ricci curvature.
In this paper, we shall also prove that in some cases,   the condition that
 $u_0$ is bounded in the above theorem can be relaxed. For example, one can
prove that if the scalar curvature has quadratic decay in the average sense,
then Theorem 3.1 is still true if the condition (i) is replaced by the condition
 that $u_0$ has sub-logarithmic growth. This is a special case of a more general
result, see Theorem 3.2.   In particular, when
$u_0(x)$ is the solution of (0.2), the Liouville result mentioned above implies
 the gap theorem proved in [C-Z] and
[N-T]. For previous results of  Liouville properties of plurisubharmonic
 functions, please see \cite{N, N-S-T1-2, N-T}.

As a by-product of our argument, we also prove a Li-Yau type differential
inequality for the positive plurisubharmonic solution $u(x,t)$ of (0.3)
  (see Theorem 2.2). Namely, we have
$$
\frac{u_t}{u}-\frac{|\nabla u|^2}{u^2}+\frac{m}{t}\ge 0
$$
exactly as in   [L-Y] for a fixed metric. Hopefully, this  differential
 inequality  will have applications to the study of the plurisubharmonic
 functions, the K\"ahler-Ricci flow and other problems.

Here is the organization of the paper. In \S1, we shall prove the
preservation
 of nonnegativity of (1,1) tensors and the linear trace LYH type inequality
for (1,1) tensors. In \S2, we shall study initial value problem (0.3) and
 prove
 the Li-Yau
type inequality for the positive solution to (0.3).
In \S3, we shall study Liouville properties of plurisubharmonic functions.

The authors would like to thank Ben Chow for his
interest and discussions and also thank Huai-Dong Cao for  discussions.

\input HP-s1
\input HP-s2
\input HP-s3

\enddocument

%% file: HP-s1.tex
\input amstex
\documentstyle{amsppt}
\magnification=1200
\hsize=13.8cm
\def\i{\sqrt{-1}}
\def\heat{\lf(\frac{\p}{\p t}-\Delta\ri)}
\def\Ric{\text{Ric}}
\def\lf{\left}
\def\ri{\right}
\def\bbar{{\bar \beta}}
\def\a{\alpha}
\def\g{\gamma}
\def\p{\partial}
\def\delbar{{\bar\delta}}
\def\ddbar{\partial\bar\partial}

\def\C{\Bbb C}

\def\Rm{{\text{\rm Rm}}}

\def \b {\beta}
\def \e{\epsilon}

\def\ba{{\bar\alpha}}
\def\bb{{\bar\beta}}
\def\abb{{\alpha\bar\beta}}

\def\tD{\widetilde \Delta}
\def\tn{\widetilde \nabla}
\def \D {\Delta}

\def\aint{\frac{\ \ }{\ \ }{\hskip -0.4cm}\int}

\subheading{\S1 A Li-Yau-Hamilton inequality}

In this section we shall prove a linear trace Li-Yau-Hamilton
 inequality which is
the K\"ahler version of
 the one obtained by Chow and Hamilton in \cite{C-H}. Just as Chow-Hamilton's
LYH inequality extends the trace Li-Yau-Hamilton inequality
 of Hamilton [H4] our differential inequality
extends the trace LYH inequality of Cao [Co1] for the scalar curvature.
Applications of this new   inequality will be given in the following sections.

Let $(M^m,g_\abb(x))$ be a complete noncompact K\"ahler manifold with bounded nonnegative holomorphic bisectional curvature. Because of the results in  \cite{Sh1}, in this section we always assume that solution of   the following Ricci-K\"ahler flow exists on $M\times[0,T]$
$$
\cases
& \frac{\p}{\p t}g_\abb=-R_\abb\\
& g_\abb(x,0)=g_\abb(x)
\endcases\tag1.1
$$
such that on $M\times[0,T]$,
\roster
\item"{(i)}" $g_\abb(x,t)$ is nonincreasing in $t$ and is  uniformly equivalent to  $g_\abb(x,0)$;
\item"{(ii)}" the curvature tensors of $g_\abb(x,t)$ are uniformly bounded;
\item"{(iii)}" there exists a constant $C$ such that
$$
|\nabla \Rm|(x,t)\le \frac{C}{t^\frac12};
$$

\item"{(iv)}" $g_\abb(x,t)$ has nonnegative holomorphic bisectional curvature.
\endroster
Sufficient conditions that (1.1) has long time existence are given in
\cite{Sh2-3, N-T}, see also \cite{C-T-Z} for the surfaces case.

In this work, we will use the maximum principle of the authors \cite{N-T, Theorem 1.2} from time to time. For the convenience of the readers, we include the statement of this maximum principle here.

\proclaim{Theorem 1.1} Let $g_{ij}(x,t)$ be a smooth family of complete Riemannian metrics defined on $M$ with $0\le t\le T$ for some $T>0$ such that     for any $T \ge t_2\ge t_1\ge0$
$$
C g_{ij}(x,t_1)\le g_{ij}(x,t_2)\le g_{ij}(x,t_1)
$$
for some   constant $C>0$ for all $x\in M$ and  let $f(x,t)$ be a smooth function such that
$(\Delta -\frac{\p}{\p t})f(x,t)\ge 0$ whenever $f(x,t)\ge0$. Assume that
  $$
\int_0^{T}\int_M \exp(-ar_0^2(x))f_{+}^2(x,s)\, dV_0\, ds <\infty
$$
for some $a>0$, where $r_0(x)$ is the distance function
to a fixed point $o\in M$ with respect to  $g_{ij}(x,0)$. Suppose $f(x,0)\le 0$ for all $x\in M$. Then $f(x,t)\le 0$ for all $(x,t)\in M\times [0,T]$.
\endproclaim

In the following let $h_{\a\bbar}(x,t)$ be a Hermitian symmetric tensor
defined on $M\times [0,T]$, which is also deformed by the complex
Lichnerowicz-Laplacian heat equation:
$$
\left(\frac{\p}{\p t}-\D \right)h_{\g\delbar}=
R_{\beta \bar{\a}\g\delbar}h_{\a\bbar}-
\frac{1}{2}\left(R_{\g\bar{p}}h_{p\delbar}+
R_{p\delbar}h_{\g\bar{p}}\right).\tag 1.2
$$

We shall obtain a LYH inequality for $h_\abb$ provided $h_\abb$ is nonnegative and does not grow very fast on $M\times[0,T]$. In application, usually we only know that $h_\abb$ is nonnegative initially. Hence we shall discuss conditions on $h_\abb$ so that nonnegativity is preserved under the flow. The following lemma is basically  from \cite{H4, Lemma 5.1}.
\proclaim{Lemma 1.1} For any $a>0$ and $C>0$ there exists a positive function $\phi(x,t)$ and $b>0$ such that $\exp(b(r_0(x)+1))\ge\phi(x,t)\ge \exp(a(r_0(x)+1))$ and
$$
\heat \phi\ge C\phi
$$
on $M\times[0,T]$, where $r_0(x)$ is the distance from a fixed point $o$ with respect to the initial metric $g(0)$.
\endproclaim
\demo{Proof} By Lemma 5.1 in \cite{H4}, there is a smooth function $f(x)$ and a constant $C_1>0$ on $M$ such that
$$
C_1^{-1}(1+r_0(x))\le f(x)\le C_1(1+r_0(x)),
$$
$$
|\nabla f|+|\nabla^2 f|\le C_1.
$$
As in \cite{H4}, we can choose $\phi(x,t)=\exp(At+\alpha f(x))$ for suitable positive constants $A$ and $\alpha$, then $\phi$ will be the required function.
\enddemo
To simplify notations, in the rest of the section, let
 $||h||$ be the norm of $h$ with respect to $g_\abb(x,t)$,
$$
\cases \Phi&=||h||^2\\
\Psi&=||\nabla h||^2=\sum_{\a\b\g}\lf(||\nabla_\g h_{\abb}||^2+||\nabla_{\bar\g}h_{\abb}||^2\ri)\\
\Lambda&=||\nabla \nabla h||^2=\sum_{\a\b\g\delta}
\lf(||\nabla_\delta\nabla_{\g} h_{\abb}||^2+||\nabla_\delta\nabla_{\bar\g}h_{\abb}||^2\ri).
\endcases\tag1.3
$$
Suppose $h$ satisfies (1.2), then direct computations show
(see \cite{H1} for example):
$$
\lf(\frac{\p}{\p t}-\D\ri)\Phi=-\Psi+A,\tag1.4
$$
$$
\lf(\frac{\p}{\p t}-\D\ri)\Psi=-\Lambda+B,\tag1.5
$$
where $A$ and $B$ satisfy the following conditions: There exists a constant
$C>0$ such that $|A|\le C\Phi$ and   $t|B|\le C\lf(\Phi+\Psi\ri)$ on
$M\times[0,T]$. Here we have used properties (ii) and (iii) of $g_\abb$.

Moreover, in   normal coordinates
$$
\split
||\nabla\Phi||^2&=\sum_\a\Phi_\a\Phi_{\bar \a}\\
&=\sum_\a\lf(\sum_{\xi,\tau}h_{\xi\bar\tau,\a}h_{\bar\xi\tau}+h_{\xi\bar\tau}h_{\bar\xi\tau,\a}\ri)\lf(\sum_{\xi,\tau}h_{\xi\bar\tau,\bar\a}h_{\bar\xi\tau}+h_{\xi\bar\tau}h_{\bar\xi\tau,\bar\a}\ri)\\
&\le 4||h||^2\sum_\a\lf[\lf(\sum_{\xi,\tau}|h_{\xi\bar\tau,\a}|^2\ri)^\frac12\lf(\sum_{\xi,\tau}|h_{\xi\bar\tau,\bar\a}|^2\ri)^\frac12\ri]\\
&\le 2||h||^2\sum_{\a,\xi\tau}\lf(|h_{\xi\bar\tau,\a}|^2+|h_{\xi\bar\tau,\bar\a}|^2\ri).
\endsplit
$$
Hence
$$
||\nabla \Phi||^2\le 2\Phi\Psi.\tag1.6
$$
Similarly,
$$
||\nabla \Psi||^2\le 2\Psi\Lambda.\tag1.7
$$
\proclaim{Lemma 1.2} Let $h_\abb$ be a tensor satisfying (1.2). Suppose
$$
||h_\abb(x,0)||\le \exp(a(1+r_0(x))\tag1.8
$$
and
$$
\int_0^T\int_M \exp(-br_0^2(x))||h||^2(x,t)dV_tdt<\infty\tag1.9
$$
for some positive constants $a$ and $b$. Then there exists a positive constant $c>0$ such that
$$
||h_\abb(x,t)||\le \exp(c(1+r_0(x))\tag1.10
$$
on $M\times [0,T]$.
\endproclaim
\demo{Proof}   By (1.4) and (1.8), it is easy to see that
$$
\heat \lf[e^{-C_1t}\lf(1+\Phi\ri)^\frac12\ri]\le 0\tag1.11
$$
for some constant $C_1>0$.

By Lemma 1.1, there exists a function $\phi(x,t)$ and constant $c>0$ such that $\exp(c(r_0(x)+1))\ge\phi(x,t)\ge \exp(a(r_0(x)+1))$ and
$$
\heat \phi\ge C\phi
$$
with $C>0$.
By (1.9), (1.10) and (1.11), we have $\phi+1\ge \lf(1+\Phi\ri)^\frac12$ by the maximum principle   Theorem 1.1. The lemma follows by choose an even larger $c$.
\enddemo

Next we shall prove that nonnegativity of $h$ will be preserved by the flow under certain conditions.
\proclaim{Proposition 1.1} Suppose $h_\abb$ satisfy (1.2) and the conditions (1.8) and (1.9) of Lemma 1.2. Suppose also that $h_\abb(x,0)\ge 0$. Then $h_\abb(x,t)\ge0$ for $t>0$.
\endproclaim
\demo{Proof} By Lemma 1.2, there exists a constant $c>0$ such that
$$
||h||(x,t)\le \exp(c(1+r_0(x))).\tag1.12
$$
By Lemma 1.1, for any $C'>0$, there exists a function $\phi$ such that
$$
\exp(c'(1+r_0(x)))\ge \phi \ge \exp(2c (1+r_0(x)))\tag1.13
$$
and
$$
\heat\phi>C'\phi\tag1.14
$$
It is enough to show that $h_{\abb}(x,t)+\epsilon\phi g_{\abb}(x,t)>0$,
for any $\epsilon >0$. Now we calculate
$$
\split
\lf(\frac{\p}{\p t}-\Delta\ri)  \lf(h_{\abb}+\epsilon\phi g_{\abb}\ri) & =
R_{\abb\g\delbar}\lf(h_{\bar{\g}\delta}
+\epsilon\phi g_{\bar{\g}\delta}\ri) -\frac{1}{2}
R_{\a \bar{p}}\lf(h_{p\bar{\beta}}+\epsilon\phi g_{p\bar{\beta}}\ri)\\
& \ \ -\frac{1}{2}R_{p\bar{\beta}}
\lf(h_{\a\bar{p}}+\epsilon\phi g_{\a\bar{p}}\ri)+\epsilon (\phi_t-\Delta \phi)
g_{\abb}-\epsilon \phi R_{\abb}.
\endsplit \tag 1.15
$$
Here we have used (1.2) and the Ricci flow equation. By (1.12), (1.13) and the fact that at $t=0$, $h_\abb+\e\phi g_\abb>0$, if  $h_{\abb}(x,t)+\epsilon\phi g_{\abb}(x,t)>0$ fails to hold at some $t>0$, then there is $(x_0,t_0)$ and unit (1,0) vector at $x_0$ with $t_0>0$ such that $\lf(h_{\abb}(x_0,t_0)+\epsilon\phi g_{\abb}(x,t)\ri)v^\a \bar v^\b=0$ and  $t_0$ is the first time that happens. As in \cite{Cw3}, we can extend $v$ in a neighborhood in space-time  of $(x_0,t_0)$ such that $\nabla v$ and $\Delta v=0$ at $(x_0,t_0)$ with respect to the metric $g(t_0)$ and such that $v$ is independent of time. Hence at $(x_0,t_0)$ we have
$$
\split
0&\ge \lf(\frac{\p}{\p t}-\Delta\ri) \lf[ \lf(h_{\abb}+\epsilon
\phi g_{\abb}\ri)v^\a \bar v^\b\ri]\\
&=\lf[\lf(\frac{\p}{\p t}-\Delta\ri)  \lf(h_{\abb}+\epsilon
\phi g_{\abb}\ri)\ri]v^\a\bar v^\b\\
&=R_{\abb\g\delbar}\lf(h_{\bar{\g}\delta}
+\epsilon\phi g_{\bar{\g}\delta}\ri)v^\a  \bar v^\b -\frac{1}{2}
R_{\a \bar{p}}\lf(h_{p\bar{\beta}}+\epsilon\phi g_{p\bar{\beta}}\ri)v^\a \bar v^\b\\
& \ \ -\frac{1}{2}R_{p\bar{\beta}}
\lf(h_{\a\bar{p}}+\epsilon\phi g_{\a\bar{p}}\ri)v^\a \bar v^\b+\epsilon (\phi_t-\Delta \phi)
g_{\abb}v^\a  \bar v^\b-\epsilon \phi R_{\abb}v^\a  \bar v^\b.
\endsplit
$$
Since $v$ minimizes $h_\abb+\e\phi g_\abb$ among all (1,0) unit vectors at $x_0$, first variation gives
$$
(h_\abb+\e\phi g_\abb)v^\a=(h_\abb+\e\phi g_\abb)\bar v^\b=0.
$$
Using also the fact that $M, g(t_0)$ has nonnegative holomorphic bisectional curvature, we conclude that
$$
0\ge \epsilon (\phi_t -\Delta \phi)-\epsilon \phi R_{\a\bar{\b}}v^a\bar v^\b>0
$$
for sufficient large $C'$, since $|Rm|$ is bounded. This   is a
contradiction.
\enddemo

We should remark that the result is still true if $M$ is compact. In this case, there is no need to impose growth condition on $h_\abb$.
Moreover if $h_\abb(x,0)$ is positive at some point, then $h_\abb(x,t)$ will be positive for all $t>0$.

In order to apply the maximum principle we also need the following estimates.
\proclaim{Lemma 1.3} Let $h_\abb$ as in Lemma 1.2. Then for any $a>0$,
$$
\int_{0}^{T}\int_M e^{-a r^2_0(x)}\Psi(x,t) dV_tdt<\infty,\tag1.16
$$
$$
\int_{0}^{T}\int_M te^{-a r^2_0(x)}\Lambda(x,t) dV_tdt<\infty\tag1.17
$$
and
$$
\int_{0}^{T}\int_M te^{-a r^2_0(x)}\Psi^2(x,t) dV_tdt<\infty.\tag1.18
$$
\endproclaim
\demo{Proof} Let $f(x)$ be a smooth function such that $0\le f\le 1$,  $f=1$ on $B_0(o,R)$, $f=0$ outside $B_0(o,2R)$ and $|\tn f|\le C/R$ for some constant $C$ independent of $R$. Here $B_0(o,R)$ is the geodesic ball with center at $o$ and radius $R$ with respect to $g(0)$. Multiply (1.4) by $f^2$ and integrating by parts, we have:
$$
\split
\int_0^{T}\int_Mf^2\Psi dV_tdt&\le -\int_0^{T}\int_Mf^2\lf(\frac{\p}{\p t}-\D\ri)\Phi dV_tdt+C_1\int_0^{T}\int_Mf^2\Phi dV_tdt\\
& \le \int_Mf^2\Phi dV_0 dt+2\int_0^{T}\int_Mf|\nabla f|\,|\nabla \Phi|dV_tdt+C_1\int_0^{T}\int_Mf^2\Phi dV_tdt
\endsplit
$$
for some constant $C_1$. Here we have used the fact that $dV_t$ is nonincreasing. Using (1.6) and Schwarz inequality, we have
$$
\int_0^{T}\int_Mf^2\Psi dV_tdt
 \le C_2\bigg[\int_Mf^2\Phi dV_0 dt+\int_0^{T}\int_M\lf(f^2+|\tn f|^2\ri)\Phi dV_0dt\bigg],
$$
for some constant $C_2$, where we have used the fact that $g(t)$ and $g(0)$ are equivalent in $[0,T]$. Hence
$$
\int_0^{T}\int_{B_0(o,R)}\Psi dV_tdt
 \le C_3\bigg[\int_{B_0(o,2R)}\Phi dV_0 dt+\int_0^{T}\int_{B_0(o,2R)}\Phi dV_0dt\bigg],\tag1.19
$$
for some constant $C_3$, where we assume $R\ge1$. By Lemma 1.1, $||h||$ is at most of exponential growth, hence it is easy to see (1.16) is true because $g(0)$ has nonnegative Ricci curvature.

To prove (1.17), multiplying (1.5) by $t f^2$ and integrating by parts we have for $R\ge 1$,
$$
\split
\int_0^T&\int_M t f^2\Lambda dV_tdt\\
&\le -\int_0^T\int_M tf^2\lf(\frac{\p}{\p t}-\D\ri)\Psi dV_tdt+C_4\int_0^T\int_M  f^2\lf(\Phi+\Psi\ri) dV_tdt\\
&\le \int_0^T\int_M  f^2\Psi dV_tdt+2\int_0^T\int_M t f|\nabla f|\,|\nabla \Psi|dV_tdt+C_4\int_0^T\int_M f^2\lf(\Phi+\Psi\ri) dV_tdt\\
&\le C_5\int_{0}^{T}\int_M\lf(f^2+|\tn f|^2\ri)\lf(\Phi+\Psi\ri)dV_tdt+\frac12\int_0^T\int_M t f^2\Lambda dV_tdt
\endsplit
$$
for some constants $C_4$, $C_5$, where we have used (1.7).
Hence if $R$ is large
$$
\int_{0}^{T}\int_{B_0(o,R)} t\Lambda dV_tdt\le 3C_5\int_{0}^{T}\int_{B_0(o,2R)}\lf(\Phi+\Psi\ri)dV_tdt.\tag1.20
$$
Combining (1.19) and (1.20), we can conclude that (1.17) is true.

To prove (1.18), multiplying (1.4) by $t f^2\Psi$ and integrating by parts, we have
$$
\split
\int_0^T&\int_M t f^2\Psi^2 dV_tdt\\
&\le -\int_0^T\int_M t f^2\Psi\lf(\frac{\p}{\p t}-\D\ri)\Phi dV_tdt+\int_0^T\int_M t f^2\Psi A dV_tdt\\
&\le C_6\int_0^T\int_M  f^2\Psi \Phi dV_tdt
+\int_0^T\int_M t  f^2\lf(\Phi\frac{\p}{\p t}\Psi+\Psi\D\Phi\ri)dV_tdt\\
&\le C_6\int_0^T\int_M  f^2\Psi\Phi  dV_tdt+\int_0^T\int_M t   f^2\lf(\Phi\D\Psi+\Psi\D\Phi\ri)dV_tdt\\
&\quad+
\int_0^T\int_M t  f^2\Phi B dV_tdt\\
&\le C_7\int_0^T\int_M  f^2\lf(\Psi+\Phi\ri)\Phi dV_tdt
+2\int_0^T\int_M t   f^2|\nabla \Phi|\,|\nabla\Psi|  dV_tdt\\
&\quad+
2\int_0^T\int_M t f|\nabla f|\lf(\Phi|\nabla \Psi|+\Psi|\nabla \Phi|\ri)dV_tdt
\endsplit
$$
for some constants $C_6$, $C_7$, where we have used (1.5). Apply (1.6) and (1.7) to $|\nabla \Phi|$ and $|\Psi|$ respectively, and use Schwarz inequality, we have,
$$
\int_0^T\int_M t f^2\Psi^2 dV_tdt\le C_8\lf(\int_0^T\int_M  \lf(f^2+|\tn f|^2\ri)\lf(\Phi+\Psi \ri)\Phi dV_tdt+tf^2\Lambda\Phi dV_tdt\ri)
$$
for some constant $C_8$. Combining this with (1.16) and (1.17) and the fact that $\Phi$   grows at most  exponentially, we can conclude that (1.18) is true.
\enddemo
\proclaim{Remark 1.1} (1.17) and (1.18) imply that for any $\e>0$,
$$
\int_{\e}^{T}\int_M te^{-a r^2_0(x)}\Lambda(x,t) dV_tdt<\infty,
$$
and
$$
\int_{\e}^{T}\int_M te^{-a r^2_0(x)}\Psi^2(x,t) dV_tdt<\infty.
$$
\endproclaim

Now we are ready to prove a LYH inequality. Let $div(h)_\a= g^{\g\delbar}\nabla_\g h_{\a\delbar}$ and $div(h)_{\bbar}=
g^{\g\delbar}\nabla_{\delbar}h_{\g \bbar}$. Consider the quantity
$$
\split
Z& =g^{\a\bbar}g^{\g\delbar}\lf[\frac{1}{2}
\lf (\nabla_{\bbar}\nabla_\g +\nabla_\g \nabla_{\bbar}\ri)
h_{\a\delbar}+R_{\a\delbar}h_{\g\bbar}+\lf(\nabla_\g h_{\a\delbar}V_{\bbar}+
\nabla_{\bbar}h_{\a\delbar}V_{\g}\ri)+h_{\a\delbar}V_{\bbar}V_\g \ri]
\\
& \ \ +\frac{H}{t}\\
&=\frac12[g^{\a\bbar}\nabla_\bbar div(h)_\a+g^{\g\delbar}\nabla_\g div(h)_\delbar]\\
&\quad+g^{\a\bbar}g^{\g\delbar}[R_{\a\delbar}h_{\g\bbar}+\nabla_\g h_{\a\delbar}V_{\bbar}+
\nabla_{\bbar}h_{\a\delbar}V_{\g}+h_{\a\delbar}V_{\bbar}V_\g ]+\frac{H}{t}\endsplit \tag 1.21
$$
where $H$ is the trace of $h_{\a\bbar}$ with respect to $g_{\a\bbar}(x,t)$.

\proclaim{Theorem 1.2} Let $h_\abb$ be a Hermitian symmetric tensor satisfying (1.2) on $M\times[0,T]$. Suppose $h_\abb(x,0)\ge0$ and satisfies (1.8) and (1.9) in Lemma 1.2.
Then $Z\ge0$ on $M\times(0,T]$ for any smooth vector field $V$ of type $(1,0)$.
\endproclaim
In order to prove the theorem, we need to compute $\lf(\frac{\p}{\p t}-\D\ri)Z$.  As in \cite{C-H}, we need to calculate $(\frac{\p}{\p t}-\Delta)Z$.  We break the computations in several lemmas.
\proclaim{Lemma 1.4} Under   normal coordinates at a point,
$$
\lf (\frac{\p }{\p t} -\Delta\ri)\lf( div(h)_\a \ri)=R_{s\bar{t}}
\nabla_t h_{\a \bar{s}}+\nabla_{\a}R_{s\bar{t}}h_{\bar{s}t}-\frac{1}{2}
R_{\a\bar{t}}\, (div(h)_t). \tag 1.22
$$
\endproclaim
\demo{Proof} Direct calculation shows
$$
\split
\frac{\p}{\p t}\lf ( g^{\g\delbar}\nabla_{\g}h_{\a\delbar}\ri) & =
\frac{\p}{\p t}\lf[g^{\g\delbar}\lf(\p_\g h_{\a\delbar}-\Gamma^{p}_{\a\g}
h_{p\delbar}\ri)\ri]\\
& = g^{\g\bar{t}}R_{s\bar{t}}g^{s\delbar}\nabla_{\g}h_{\a\delbar} +
g^{\g\delbar}\nabla_{\g}(\frac{\p}{\p t}h_{\a\delbar})-g^{\g\delbar}
\lf(\frac{\p }{\p t}\Gamma^p_{\a\g}\ri) h_{p\delbar}\\
& = R_{s\bar{t}}\nabla_{t}h_{\a\bar{s}}+
\nabla_{\g}\lf(\D h_{\a\bar{\g}}+R_{\a\bar{\g} s\bar{t}}h_{\bar{s}t}
-\frac{1}{2}R_{\a\bar{t}}h_{t\bar{\g}}-
\frac{1}{2}R_{t\bar{\g}}h_{\a\bar{t}}\ri)\\
& \ \ +\nabla_{\g}R_{\a\bar{p}}
h_{p\bar{\g}}\\
&= R_{s\bar{t}}\nabla_{t}h_{\a\bar{s}}+\nabla_\a R_{s\bar{t}}h_{\bar{s}t}
+R_{\a\bar{\g}s\bar{t}}\nabla_\g h_{\bar{s}t}+
\frac{1}{2}\nabla_{\g}R_{\a\bar{t}}h_{t\bar{\g}}\\
& \ \ -\frac{1}{2}R_{\a\bar{t}}\nabla_{\g}h_{t\bar{\g}}-
\frac{1}{2}\nabla_t Rh_{\a\bar{t}}-
\frac{1}{2}R_{t\bar{\g}}\nabla_{\g}h_{\a\bar{t}}+\nabla_{\g}
(\D h_{\a\bar{\g}}).
\endsplit \tag1.23
$$
Now we calculate $\nabla_{\g}(\D h_{\a\bar{\g}})$. By definition,
$$
\nabla_{\g}(\D h_{\a\bar{\g}}) =\frac{1}{2}\nabla_{\g}
\lf(\nabla_s\nabla_{\bar{s}}+\nabla_{\bar{s}}\nabla_s\ri)h_{\a\bar{\g}}.
$$
On the other hand,
$$
\split
\nabla_{\g} \nabla_{s} \nabla_{\bar{s}} h_{\a\bar{\g}} & =
\nabla_s \nabla_\g \nabla_{\bar{s}} h_{\a\bar{\g}} \\
& = \nabla_{s} \lf[ \nabla_{\bar{s}} \nabla_\g  h_{\a\bar{\g}}-
R_{\a\bar{p}\g\bar{s}}h_{p\bar{\g}}+R_{p\bar{\g}\g\bar{s}} h_{\a\bar{p}}\ri]\\
&= \nabla_{s}\nabla_{\bar{s}} \nabla_\g h_{\a\bar{\g}}-
\nabla_\g R_{\a\bar{p}}h_{p\bar{\g}}-R_{\a\bar{p}\g\bar{s}}
\nabla_s h_{p\bar{\g}}+\nabla_p R h_{\a\bar{p}}+
R_{p\bar{s}}\nabla_s h_{\a\bar{p}}.
\endsplit
$$
Similarly,
$$
\split
\nabla_\g \nabla_{\bar{s}}\nabla_s h_{\a\bar{\g}} & = \nabla_{\bar{s}}
\nabla_\g\nabla_s h_{\a\bar{\g}}+R_{p\bar{\g}\g\bar{s}}
\nabla_s h_{\a\bar{p}}-R_{s\bar{p}\g\bar{s}}\nabla_p h_{\a\bar{\g}}\\
& \ \ \ - R_{\a\bar{p}\g \bar{s}}\nabla_s h_{p\bar{\g}}\\
&= \nabla_{\bar{s}}\nabla_s\nabla_{\g}h_{\a\bar{\g}}+
R_{p\bar{s}}\nabla_s h_{\a\bar{p}}-R_{\g\bar{p}}\nabla_p h_{\a\bar{\g}}
-R_{\a\bar{p}\g \bar{s}}\nabla_s h_{p\bar{\g}}.
\endsplit
$$
Combining the above three we have
$$
\split
\nabla_{\g}(\D h_{\a\bar{\g}}) &= \D(\nabla_\g h_{\a\bar{\g}})-\frac{1}{2}
\nabla_\g R_{\a\bar{p}}h_{p\bar{\g}}-
R_{\a\bar{p}\g\bar{s}}\nabla_s h_{p\bar{\g}} \\
& \ \ \ +\frac{1}{2}\nabla_p R h_{\a\bar{p}}+
R_{p\bar{s}}\nabla_s h_{\a\bar{p}}-\frac{1}{2}R_{\g\bar{p}}\nabla_p
h_{\a\bar{\g}}.
\endsplit
$$
Plugging the above into (1.23), the lemma is proved.
\enddemo

\proclaim{Lemma 1.5} Under   normal coordinates at a point,
$$
\split
\lf(\frac{\p}{\p t}-\D\ri)\lf(g^{\abb}\nabla_{\bar{\b}} div(h)_\a\ri)& =
R_{s\bar{\a}}\nabla_{\bar{s}}div(h)_\a+
\nabla_{\bar{\a}}R_{s\bar{t}}\nabla_th_{\a\bar{s}}+\nabla_{\a}
R_{s\bar{t}}\nabla_{\bar{\a}}h_{\bar{s}t}\\
&\ \ \ + R_{s\bar{t}}\nabla_{\bar{\a}}\nabla_t h_{\a\bar{s}}
+\lf(\nabla_{\bar{\a}}\nabla_{\a}R_{s\bar{t}}\ri) h_{\bar{s}t}.
\endsplit \tag1.24
$$
\endproclaim
\demo{Proof} Direct calculation shows that
$$
\split
\frac{\p}{\p t}\lf(g^{\a\bbar}\nabla_{\bbar}div(h)_\a\ri) & = g^{\a\bar{t}}
R_{s\bar{t}}g^{s\bbar}\nabla_{\bbar}div(h)_\a +
\nabla_{\bar{\a}}\lf(\frac{\p}{\p t} div(h)_\a\ri)\\
& = R_{s\bar{\a}}\nabla_{\bar{s}}\lf(div(h)_\a\ri)\\
& \ \ + \nabla_{\bar{\a}}\lf[ \D div(h)_\a +R_{s\bar{t}}
\nabla_t h_{\a\bar{s}}+\nabla_{\a}R_{s\bar{t}}h_{\bar{s}t}-
\frac{1}{2}R_{\a\bar{t}}div(h)_t\ri],
\endsplit
$$
by Lemma 1.4. Therefore we have that
$$
\split
\frac{\p}{\p t}\lf(g^{\a\bbar}\nabla_{\bbar}div(h)_\a\ri) &=
\nabla_{\bar{\a}}\lf(\D div(h)_\a\ri)+\frac{1}{2}R_{s\bar{\a}}
\nabla_{\bar{s}}\lf(div(h)_\a\ri)+\nabla_{\bar{\a}}R_{s\bar{t}}\nabla_t
h_{\a\bar{s}}\\
& \ \ +\nabla_{\a}R_{s\bar{t}}\nabla_{\bar{\a}}h_{\bar{s}t} +
 R_{s\bar{t}}\nabla_{\bar{\a}}\nabla_t h_{\a\bar{s}}
+\nabla_{\bar{\a}}\nabla_\a R_{s\bar{t}}h_{\bar{s}t}-
\frac{1}{2}\nabla_{\bar{t}}R\lf( div(h)_t\ri).
\endsplit \tag1.25
$$
Now we calculate $\nabla_{\bar{\a}}\lf(\D div(h)_\a\ri)$.
By definition
$$
\nabla_{\bar{\a}}\lf(\D div(h)_\a\ri)=\frac{1}{2}\nabla_{\bar{\a}}
\nabla_s\nabla_{\bar{s}}div(h)_{\a}+\frac{1}{2}\nabla_{\bar{\a}}
\nabla_{\bar{s}}\nabla_s div(h)_\a.
$$
On the other hand
$$
\split
\nabla_{\bar{\a}}\nabla_{\bar{s}}\nabla_s div(h)_\a &=
\nabla_{\bar{s}}\nabla_{\bar{\a}}\nabla_s div(h)_\a \\
& = \nabla_{\bar{s}}\lf[ \nabla_s\nabla_{\bar{\a}}div(h)_\a+
R_{\a\bar{p}s\bar{\a}}div(h)_p\ri]\\
&= \nabla_{\bar{s}}\nabla_s\nabla_{\bar{\a}}div(h)_\a+
\lf(\nabla_{\bar{s}}R\ri)\lf(
div(h)_s\ri)+R_{s\bar{p}}\nabla_{\bar{s}}div(h)_p
\endsplit
$$
and
$$
\split
\nabla_{\bar{\a}}\nabla_s\nabla_{\bar{s}}div(h)_\a &=
\nabla_{s}\nabla_{\bar{\a}}\nabla_{\bar{s}}div(h)_{\a} +
R_{\bar{p}s}\nabla_{\bar{s}} div(h)_p-R_{p\bar{\a}}\nabla_{\bar{p}}div(h)_\a\\
&= \nabla_{s}\nabla_{\bar{s}}\nabla_{\bar{\a}}div(h)_{\a}+
R_{\bar{p}s}\nabla_{\bar{s}} div(h)_p-R_{p\bar{\a}}\nabla_{\bar{p}}div(h)_\a.
\endsplit
$$
Combining the above three  we have that
$$
\nabla_{\bar{\a}}\lf(\D div(h)_\a\ri)=\D\lf(\nabla_{\bar{\a}}div(h)_{\a}\ri)
+\frac{1}{2}\nabla_{\bar{s}}R(div(h)_s)+\frac{1}{2}R_{s\bar{p}}
\nabla_{\bar{s}}div(h)_p.
$$
Plugging into (1.25), this completes the proof of Lemma 1.3.
\enddemo

Taking the conjugation we will have the following lemmas.

\proclaim{Lemma 1.4\'}
Under   normal coordinates at a point,
$$
\lf(\frac{\p}{\p t}-\D\ri)\lf(div (h)_{\bbar}\ri)= R_{s\bar{p}}
\nabla_{\bar{s}}h_{p\bbar}+\nabla_{\bbar}R_{p\bar{\g}}h_{\g\bar{p}}
-\frac{1}{2}R_{t\bbar}div(h)_{\bar{t}} \tag1.26
$$
\endproclaim
\proclaim{Lemma 1.5\'} Under  normal coordinates at a point,
$$
\split
\lf(\frac{\p}{\p t}-\D\ri)\lf(g^{\b\bar\a}\nabla_\b div(h)_{\bar{\a}}\ri)& =
R_{\a\bar{p}}\nabla_p div(h)_{\bar{\a}}+\nabla_{\a}R_{s\bar{p}}
\nabla_{\bar{s}}h_{p\bar{\a}}
+\nabla_{\bar{\a}}R_{p\bar{\g}}\nabla_\a h_{\g \bar{p}}\\
& \ \ \ + R_{s\bar{p}}\nabla_{\a}\nabla_{\bar{s}}h_{p\bar{\a}}+
\lf(\nabla_{\a}\nabla_{\bar{\a}}R_{p\bar{\g}}\ri)h_{\g\bar{p}}
\endsplit
\tag1.27
$$
\endproclaim

Now we are ready to calculate $\lf(\frac{\p}{\p t}-\D\ri)Z$. By Proposition 1.1, $h$ is   nonnegative. However $h$ may be zero somewhere,  we consider $\widehat Z$ instead, where :
$$
\widehat Z=\text{\rm I}+\text{\rm II}+\text{\rm III}+\text{\rm IV}+\text{\rm V},
$$
where
$$
\split
& \text{\rm I}  =\frac12[g^{\a\bbar}\nabla_\bbar div(h)_\a)+g^{\g\delbar}\nabla_\g div(h)_\delbar], \\
& \text{\rm II}  = g^{\a\bbar}g^{\g\delbar}[R_{\a\delbar}h_{\g\bbar}+\epsilon \Cal R],\\
& \text{\rm III}  = g^{\a\bbar}div(h)_\a V_{\bar{\b}}+
g^{\g\delbar}div(h)_{\bar{\delta}}V_\g,\\
& \text{\rm IV} = g^{\a\bbar}g^{\g\delbar}\widetilde h_{\a\bar\delta}V_{\bar{\b}}V_{\a},\\
& \text{\rm V}  =\frac{H+\e m}{t}.
\endsplit
$$
where $\e>0$ is a fixed constant, $\Cal R$ is the scalar curvature, and $\widetilde h_\abb=h_\abb+\e g_\abb$.
We calculate them one by one. In the following, we always do computations in a normal coordinates at a point because the final result will not depend on the choice of coordinates.

From Lemma 1.5, Lemma 1.5\' \ and the second
Bianchi identity we have that
$$
\split
\lf(\frac{\p}{\p t}-\D \ri)\text{\rm  I} &= \frac{1}{2}\lf[ R_{\a\bar{p}}
\nabla_{p} div(h)_{\bar{\a}}+R_{p\bar{\a}}\nabla_{\bar{p}}div(h)_{\a} \ri] +
\frac{1}{2}\lf[ R_{s\bar{p}}\nabla_\a\nabla_{\bar{s}}h_{p\bar{\a}}+R_{s\bar{p}}\nabla_{\bar{\a}}\nabla_{p}h_{\a\bar{s}} \ri]  \\
&\ \ \ +\D R_{s\bar{t}}h_{\bar{s}t}+
\nabla_{\bar{t}}R_{s\bar{\a}}\nabla_{t}h_{\a \bar{s}}+
\nabla_t R_{s\bar{\a}}\nabla_{\bar{t}}h_{\bar{s}\a}.
\endsplit \tag1.28
$$

$$
\split
\lf(\frac{\p}{\p t}-\D \ri) \text{\rm II} &= \lf(\frac{\p}{\p t}-\D \ri)\lf(
g^{\a\delbar}g^{\g\bbar}R_{\a\bbar}h_{\g\delbar}+\e \Cal R\ri)\\
& = 2R_{\beta \bar{\g}}R_{\a\bbar}h_{\g\bar{\a}}+
\lf(\lf(\frac{\p}{\p t}-\D \ri)R_{\a\bbar}\ri)h_{\beta\bar{\a}} +R_{\a\bbar}\lf(\lf(\frac{\p}{\p t}-\D \ri)h_{\beta\bar{\a}} \ri)\\
& \ \  \-
\nabla_{s}R_{\a\bbar}\nabla_{\bar{s}}h_{\beta\bar{\a}}-
\nabla_{\bar{s}}R_{\a\bbar}\nabla_s h_{\beta\bar{\a}}+\e|R_\abb|^2\\
& = 2R_{\a\bbar s\bar{t}}R_{\bar{s}t}h_{\beta\bar{\a}}-
\nabla_s R_{\a\bbar}\nabla_{\bar{s}}h_{\beta\bar{\a}}
-\nabla_{\bar{s}}R_{\a\bbar}\nabla_s h_{\beta\bar{\a}}+\e|R_\abb|^2.
\endsplit \tag1.29
$$
Here we have used (1.2) and the equation satisfies by the Ricci form \cite{Sh3}:
$$
\lf(\frac{\p}{\p t}-\D \ri)R_{\a\bbar}=R_{\a\bbar\g\delbar}
R_{\bar{\g}\delta}-R_{\a\bar{s}}R_{s\bbar}. \tag1.30
$$
Using Lemma 1.4, Lemma 1.4\' \ and the second Bianchi identity we have that
$$
\split
\lf(\frac{\p}{\p t}-\D \ri)\text{\rm III} & = \lf(\frac{\p}{\p t}-\D \ri)\lf[
g^{\a\bbar}\lf(div(h)_\a V_{\bbar} +div(h)_{\bbar}V_\a \ri)\ri] \\
& = R_{\bar{\a}\beta}div(h)_\a V_{\bbar}+
R_{\bar{\a}\beta}div(h)_{\bbar}V_\a \\
&\ \ \ + div(h)_\a \lf(
\lf(\frac{\p}{\p t}-\D \ri)V_{\bar{\a}}\ri)+div(h)_{\bar{\a}}
\lf(\lf(\frac{\p}{\p t}-\D \ri)V_\a\ri)\\
&\ \ \
+ \lf( R_{s\bar{t}}\nabla_{t}h_{\a\bar{s}}+
\nabla_{s}R_{\a\bar{t}}h_{\bar{s}t}-\frac{1}{2}R_{\a\bar{t}}div(h)_t\ri)
V_{\bar{\a}}\\
&\ \ \
+\lf(R_{s\bar{t}}\nabla_{\bar{s}}h_{t\bar{\a}}
+\nabla_{\bar{s}}R_{t\bar{\a}}h_{s\bar{t}}-\frac{1}{2}R_{t\bar{\a}}
div(h)_{\bar{t}}\ri)V_\a \\
& \ \ \ -\nabla_s div(h)_\a\nabla_{\bar{s}}V_{\bar{\a}}-\nabla_{\bar{s}}
div(h)_\a\nabla_s V_{\bar{\a}}\\
&\ \ \ -\nabla_s div(h)_{\bar{\a}}
\nabla_{\bar{s}}V_{\a}-\nabla_{\bar{s}}div(h)_{\bar{\a}}\nabla_s V_\a.
\endsplit \tag1.31
$$
Using (1.2) we have
$$
\split
\lf(\frac{\p}{\p t}-\D \ri)\text{\rm IV} &= \lf(\frac{\p}{\p t}-\D \ri)
\lf(g^{\a\bbar}g^{\g\delbar}\widetilde h_{\a\delbar}V_\g V_{\bbar}\ri) \\
&= R_{\a\bbar s\bar{t}}h_{\bar{s}t}V_{\beta}V_{\bar{\a}}+
\frac{1}{2}R_{\a\bar{s}}h_{s\bar{\g}}V_{\g}V_{\bar{\a}}
+\frac{1}{2}h_{\a\bar{s}}R_{s\bar{\g}}V_\g V_{\bar{\a}}\\
& \ \ \ +\widetilde h_{\a\bar{\g}}\lf(\lf(\frac{\p}{\p t}-\D \ri)V_\g\ri)V_{\bar{\a}}
+\widetilde h_{\a\bar{\g}}V_\g \lf(\lf(\frac{\p}{\p t}-\D \ri)V_{\bar{\a}}\ri)\\
& \ \ \
-\nabla_s h_{\a\bar{\g}}\nabla_{\bar{s}}\lf(V_\g V_{\bar{\a}}\ri)
-\nabla_{\bar{s}}h_{\a\bar{\g}}\nabla_{s}\lf(V_{\g}V_{\bar{\a}}\ri)\\
& \ \ \ - \widetilde h_{\a\bar{\g}}\lf[\nabla_s V_\g \nabla_{\bar{s}}V_{\bar{\a}}+
\nabla_{\bar{s}}V_{\g}\nabla_{s}V_{\bar{\a}}+\e R_{\a\bar\g}V_{\bar\a}V_\g\ri].
\endsplit
\tag1.32
$$
Taking trace on (1.2) one can have
$$
\lf(\frac{\p}{\p t}-\D \ri)\text{\rm V}=\frac{R_{\a\bar{s}}h_{s\bar{\a}}}{t}-
\frac{H+\e m}{t^2}. \tag1.33
$$
Now combining them together we have that
$$
\split
\lf(\frac{\p}{\p t}-\D \ri)Z & = Y_1+\frac{1}{2}
\lf[ R_{\a\bar{p}}\nabla_p div(h)_{\bar{\a}}+R_{p\bar{\a}}
\nabla_{\bar{p}}div(h)_{\a} \ri] \\
& \ \ \ +\frac{1}{2}\lf[R_{s\bar{p}}\nabla_{\a}\nabla_{\bar{s}}
h_{p\bar{\a}}+R_{s\bar{p}}\nabla_{\bar{\a}}\nabla_{p}h_{\a\bar{s}}\ri]\\
& \ \ \ +
 R_{\a\bbar s \bar{t}}R_{\bar{\a}\beta}h_{t\bar{s}}
+R_{\bar{\a}\beta}div(h)_\a V_{\bbar}+R_{\bar{\a}\beta}div(h)_{\bbar}V_\a \\
& \ \ \  + div(h)_\a \lf(
\lf(\frac{\p}{\p t}-\D \ri)V_{\bar{\a}}\ri)+div(h)_{\bar{\a}}
\lf(\lf(\frac{\p}{\p t}-\D \ri)V_\a\ri)\\
& \ \ \ -\nabla_s div(h)_\a\nabla_{\bar{s}}V_{\bar{\a}}-\nabla_{\bar{s}}
div(h)_\a\nabla_s V_{\bar{\a}}\\
&\ \ \ -\nabla_s div(h)_{\bar{\a}}
\nabla_{\bar{s}}V_{\a}-\nabla_{\bar{s}}div(h)_{\bar{\a}}\nabla_s V_\a\\
& \ \ \ + R_{s\bar{t}}\nabla_t h_{\a\bar{s}}V_{\bar{\a}}-
\frac{1}{2}R_{\a\bar{t}}div(h)_t V_{\bar{\a}}+R_{s\bar{t}}\nabla_{\bar{s}}
h_{t\bar{\a}}V_\a  -\frac{1}{2}R_{t\bar{\a}}div(h)_{\bar{t}}V_\a
\\
& \ \ \ + \frac{1}{2}R_{\a\bar{s}}h_{s\bar{\g}}V_{\g}V_{\bar{\a}}
+\frac{1}{2}h_{\a\bar{s}}R_{s\bar{\g}}V_\g V_{\bar{\a}}\\
& \ \ \  + \widetilde h_{\a\bar{\g}}\lf(\lf(\frac{\p}{\p t}-\D \ri)V_\g\ri)V_{\bar{\a}}
+\widetilde h_{\a\bar{\g}}V_\g \lf(\lf(\frac{\p}{\p t}-\D \ri)V_{\bar{\a}}\ri)\\
& \ \ \
-\nabla_s h_{\a\bar{\g}}\nabla_{\bar{s}}\lf(V_\g V_{\bar{\a}}\ri)
-\nabla_{\bar{s}}h_{\a\bar{\g}}\nabla_{s}\lf(V_{\g}V_{\bar{\a}}\ri)\\
& \ \ \ - \widetilde h_{\a\bar{\g}}\lf[\nabla_s V_\g \nabla_{\bar{s}}V_{\bar{\a}}+
\nabla_{\bar{s}}V_{\g}\nabla_{s}V_{\bar{\a}}\ri] -\frac{H+\e m}{t^2}+\e |R_\abb|^2+\e R_{\a\bar\g}V_{\bar\a}V_\g,
\endsplit
\tag1.34
$$
where
$$
Y_1=\lf[ \D R_{s\bar{t}} +R_{s\bar{t}\a\bbar}R_{\bar{\a}\beta}
+\nabla_\a R_{s\bar{t}}V_{\bar{\a}}+\nabla_{\bar{\a}}R_{s\bar{t}}V_{\a}+
R_{s\bar{t}\a\bbar}V_{\bar{\a}}V_\beta +\frac{R_{s\bar{t}}}{t}\ri]h_{\bar{s}t}.
\tag1.35
$$
By Proposition 1.1, $h_\abb\ge0$ on $M\times[0,T]$. Hence by Cao's
LYH inequality \cite{Co1-2} and the fact that  $(M, g_{\a\bbar}(x,t))$ has   nonnegative holomorphic bisectional
curvature, the two factors in $Y_1$ are all nonnegative tensors. Therefore
$Y_1\ge 0$. Since  $\widetilde h_{\a\bbar}\ge \e g_\abb$, for each $(x,t)$, $\widehat Z$ attains minimum for some $V$. Then by the first variation we have
$$
div(h)_\a +\widetilde h_{\a\bar{\g}}V_\g =0\  \text{  and  }\  div(h)_{\bar{\a}}
+\widetilde h_{\g\bar{\a}}V_{\bar{\g}}=0. \tag1.36
$$
Direct calculation also shows that
$$
\split
R_{p\bar{\a}}\nabla_{\bar{p}}div(h)_\a+R_{\a\bar{p}}\nabla_p div(h)_{\bar{\a}} &=
R_{s\bar{p}}\nabla_{\a}\nabla_{\bar{s}}h_{p\bar{\a}}+
R_{s\bar{p}}\nabla_{\bar{\a}}\nabla_{p}h_{\a\bar{s}}+2R_{\a\bar{p}}
R_{p\bar{\a}s\bar{\g}}h_{\g\bar{s}}\\
& \ \ \  -2R_{\a\bar{p}}R_{p\bar{s}}h_{s\bar{\a}} \endsplit\tag1.37
$$
Combining (1.34)-(1.37) we have that
$$
\split
\lf(\frac{\p}{\p t}-\D \ri)\widehat Z & = Y_1+
\lf[ R_{\a\bar{p}}\nabla_p div(h)_{\bar{\a}}+R_{p\bar{\a}}
\nabla_{\bar{p}}div(h)_{\a} \ri]+R_{\a\bar{p}}R_{p\bar{s}}h_{s\bar{\a}} \\
& \ \ \ -\nabla_s div(h)_\a\nabla_{\bar{s}}V_{\bar{\a}}-\nabla_{\bar{s}}
div(h)_\a\nabla_s V_{\bar{\a}}\\
&\ \ \ -\nabla_s div(h)_{\bar{\a}}
\nabla_{\bar{s}}V_{\a}-\nabla_{\bar{s}}div(h)_{\bar{\a}}\nabla_s V_\a\\
& \ \ \ +R_{s\bar{t}}\nabla_t h_{\a\bar{s}}V_{\bar{\a}}+
R_{s\bar{t}}\nabla_{\bar{s}}
h_{t\bar{\a}}V_\a
-\nabla_s h_{\a\bar{\g}}\nabla_{\bar{s}}\lf(V_\g V_{\bar{\a}}\ri)
-\nabla_{\bar{s}}h_{\a\bar{\g}}\nabla_{s}\lf(V_{\g}V_{\bar{\a}}\ri)\\
& \ \ \ - \widetilde h_{\a\bar{\g}}\lf[\nabla_s V_\g \nabla_{\bar{s}}V_{\bar{\a}}+
\nabla_{\bar{s}}V_{\g}\nabla_{s}V_{\bar{\a}}\ri] -\frac{H+\e m}{t^2}
\endsplit
\tag1.38
$$
Differentiate (1.36) we have
$$
\split
& \nabla_s div(h)_\a+\lf(\nabla_s h_{\a\bar{\g}}\ri)V_{\g}
+\widetilde h_{\a\bar{\g}}\nabla_s V_\g =0, \ \ \nabla_s div(h)_{\bar{\a}}
+\lf(\nabla_s h_{\g\bar{\a}}\ri) V_{\bar{\g}}+\widetilde h_{\g\bar{\a}}\nabla_s
V_{\bar{\g}}=0,\\
& \nabla_{\bar{s}}div(h)_\a+\lf(\nabla_{\bar{s}}h_{\a\bar{\g}}\ri)V_{\g}
+\widetilde h_{\a\bar{\g}}\nabla_{\bar{s}} V_\g =0,\ \ \nabla_{\bar{s}} div(h)_{\bar{\a}}
+\lf(\nabla_{\bar{s}} h_{\g\bar{\a}}\ri) V_{\bar{\g}}+
\widetilde h_{\g\bar{\a}}\nabla_{\bar{s}} V_{\bar{\g}}=0.
\endsplit \tag1.39
$$
Plugging this above one into (1.34) we have that
$$
\split
\lf(\frac{\p}{\p t}-\D \ri)\widehat Z & = Y_1+R_{\a\bar{p}}R_{p\bar{s}}h_{s\bar{\a}}
-R_{\a\bar{p}}h_{\g\bar{\a}}\nabla_{p}V_{\bar{\g}}-R_{p\bar{\a}}h_{\a\bar{\g}}
\nabla_{\bar{p}}V_{\g}\\
&\ \ \ +\widetilde h_{\g\bar{\a}}\nabla_s V_{\bar{\g}}\nabla_{\bar{s}}V_{\a}
+\widetilde h_{\g\bar{\a}}\nabla_{\bar{s}}V_{\bar{\g}}\nabla_s V_\a -\frac{H+\e m}{t^2}+\e|R_\abb|^2.
\endsplit
\tag1.40
$$
Let
$$
Y_2=\widetilde h_{\g\bar{\a}}\lf[ \nabla_p V_{\bar{\g}}-R_{p\bar{\g}}
-\frac{1}{t}g_{p\bar{\g}}\ri]\lf[\nabla_{\bar{p}}V_\a-R_{\a\bar{p}}
-\frac{1}{t}
g_{\bar{p}\a}\ri]+\widetilde h_{\g\bar{\a}}\nabla_{\bar{p}}V_{\bar{\g}}\nabla_pV_{\a}.
\tag1.41
$$
By Proposition 1.1 again, $Y_2\ge0$.
$$
\lf(\frac{\p}{\p t}-\D \ri)\widehat Z  = Y_1+Y_2-\frac{1}{t}
\lf[-\widetilde h_{\g\bar{\a}}\nabla_\a V_{\bar{\g}}-\widetilde h_{\g\bar{\a}}\nabla_{\bar{\a}}V_{\g}
+2R_{\a\bar{\g}}h_{\g\bar{\a}}+\frac{2(H+\e m)}{t}+2\e \Cal R\ri]. \tag1.42
$$
Using (1.36) we also know that
$$
\widehat Z=R_{\a\bbar}h_{\bar{\a}\beta}-\frac{1}{2}\widetilde h_{\a\bbar}\nabla_{\bar{\a}}V_{\beta}
-\frac{1}{2}\widetilde h_{\beta\bar{\a}}\nabla_{\a}V_{\bbar}+\frac{H+\e m}{t}+\e\Cal R. \tag 1.43
$$
Plugging into (1.42) and using the fact that $Y_1\ge 0$ and $Y_2\ge0$, we have
$$
\lf(\frac{\p}{\p t}-\D \ri)\widehat Z \ge -\frac{2\widehat Z}{t}. \tag1.44
$$
where $V$ is the smooth vector field given by (1.36). Note that both sides of (1.44) do not depend on the choice of coordinates.

\demo{Proof of Theorem 1.2} Since $\widetilde h_\abb\ge \epsilon g_\abb$ on $M\times[0,T]$, by (1.36) and (1.39), we have
$$
||V||\le C_1||\nabla h||,
$$
and
$$
||\nabla V||\le C_2\lf(||\nabla\nabla h||+||\nabla h||^2\ri),
$$
for some constants $C_1$ and $C_2$. Combining this with (1.44),
we have
$$
|t^2\widehat Z|^2\le C_3\lf(\Phi +\Phi(\Psi^2+\Lambda)+1\ri)\tag1.45
$$
for some constant $C_3$.   By (1.43), the corresponding $\widehat Z$
 satisfies  $$
\lf(\frac{\p}{\p t}-\D \ri)(t^2\widehat Z)\ge0 \tag1.46
$$
for the vector field which minimizes $\widehat Z$. By  Lemma 1.2, Lemma 1.3 and (1.45), we have
$$
\int_0^{T}\int_M\exp(-ar^2_0(x))\lf(t^2\widehat Z\ri)^2dV_tdt<\infty
$$
for any $a>0$. By the maximum principle Theorem 1.1, we have $t^2\widehat Z\ge0$ because it is obvious that $t^2\widetilde Z=0$ at $t=0$. Since this is true for
the vector field $V$ minimizing $\widehat Z$, we have $\widehat Z\ge0$ for any
 (1,0) vector field. Let $\epsilon\to0$ and the proof of the theorem is completed.
\enddemo

\noindent
{\bf Remark 1.1.} (i) The theorem is still true for the case that $M$ is compact with positive holomorphic bisectional curvature because of the result in \cite{Co1}. (ii) When $h_{\abb}=R_{\abb}$, it is known that the Ricci tensor
satisfies (1.2). Therefore we can apply Theorem 1.1 to this case. Since
$$
div(h)_{\a}=\nabla_{\g}R_{\a\bar{\g}}=\nabla_{\a}{\Cal R} \ \ \text{ and }
\ \ div(h)_{\delbar}= \nabla_{\bar{\a}}R_{\a\delbar}=\nabla_{\delbar}{\Cal R}
$$
we have
$$
Z=\D {\Cal R} +R_{\a\bbar}R_{\bar{\a}\beta}+\nabla_{\a}{\Cal R}
V_{\bar{\a}}+\nabla_{\bar{\a}}{\Cal R}V_{\a}+R_{\a\bbar}V_{\bar{\a}}V_{\beta}
+\frac{{\Cal R}}{t}\ge 0. \tag 1.46
$$
It is the trace of the LYH inequality
 proved by Cao in [Co1-2]. Hence Theorem 1.1 can be considered as a generalization of the LYH inequality of Cao for the scalar curvature.
However, we should emphasis that   Cao's result has been used in the proof of  Theorem 1.2.


%% file: HP-s2.tex
\input amstex
\documentstyle{amsppt}
\magnification=1200
\hsize=13.8cm
\def\i{\sqrt{-1}}
\def\Ric{\text{Ric}}
\def\lf{\left}
\def\ri{\right}
\def\bbar{{\bar \beta}}
\def\a{\alpha}
\def\g{\gamma}
\def\p{\partial}
\def\delbar{{\bar\delta}}
\def\ddbar{\partial\bar\partial}

\def\C{\Bbb C}

\def\Rm{{\text{\rm Rm}}}

\def \b {\beta}
\def \e{\epsilon}

\def\ba{{\bar\alpha}}
\def\bb{{\bar\beta}}
\def\abb{{\alpha\bar\beta}}

\def\tD{\widetilde \Delta}
\def\tn{\widetilde \nabla}
\def \D {\Delta}

\def\heat{\lf(\frac{\p}{\p t}-\Delta\ri)}
\def\aint{\frac{\ \ }{\ \ }{\hskip -0.4cm}\int}
\subheading{\S2 Deforming plurisubharmonic functions}

Let $(M^m,g_\abb(x,t))$ be a complete noncompact K\"ahler manifold with
bounded nonnegative holomorphic bisectional curvature deformed by the
K\"ahler-Ricci flow (1.1). As in last section we assume that
(1.1) has solution on $M\times [0,T]$ which satisfies conditions (i)-(iv) in that section.
In this section we shall study the plurisubharmonic functions deformed
by the time-dependent heat equation:
$$
\cases
 &\lf(\frac{\p}{\p t}-\Delta\ri)u(x,t) =0,\\
&\  u(x,t) =u_0(x)
\endcases \tag 2.1
$$
where $\Delta =g^{\abb}(x,t)\frac{\p^2}{\p z_{\a}\p \bar{z}_{\beta}}$ and
$u_0(x)$ is a plurisubharmonic function on $M$.

First, we shall consider more general case and drop the assumption
that $u_0$ is plurisubharmonic. We have the following existence result.

\proclaim{Proposition 2.1} Let $u_0$ be a continuous function such that $|u_0(x)|\le  \exp(a\lf(r_0(x)+1)\ri)$ for all $x$ for some positive constant $a>0$. Then
there is a unique solution of (2.1) on $M\times [0,T]$ such that
$$
|u(x,t)|\le  \exp(b(r_0(x)+1))\tag2.2
$$
on $M\times[0,T]$ for some positive constant $b$.
\endproclaim
\demo{Proof} By Lemma 1.1, there exists a function $\varphi(x)$ such that
$$ \exp(b(r_0(x)+1))\ge\varphi(x)\ge  \exp(a(r_0(x)+1))
$$
for some positive constant  and $b$ for all $(x,t)\in M\times [0,T]$, and
$$
 \heat \varphi\ge0.
$$
Using $\varphi$ and $-\varphi$ as barriers, the existence part of the proposition follows. Uniqueness follows from the maximum principle Theorem 1.1.
\enddemo
Next we shall study   properties of the solution $u$ obtained in the proposition. We need the following lemma.
\proclaim{Lemma 2.1}
Let $u(x,t)$ be a solution of (2.1).
Then $u_{\abb}$ satisfies the complex Lichnerowicz heat equation (1.2).
\endproclaim
\demo{Proof}   Differentiate (2.1) we have
$$
(u_t)_{\g\delbar}=R_{\beta \bar{\a}\g\delbar}u_{\a\bbar}+
g^{\a\bbar}u_{\a\bbar\g\delbar}. \tag 2.3
$$
By definition $\Delta u_{\a\bbar}=
\frac{1}{2}\left(u_{\a\bbar,\g\bar{\g}}+u_{\a\bbar, \bar{\g}\g}\right)$,
in  normal coordinates at a point.
We  need to calculate the difference between the partial derivative
$u_{\a\bbar\g\delbar}$ and the
covariant derivative
$u_{\a\bbar,\g\delbar}$.
Direct computations show that, for normal coordinates at a point
$$
u_{\g\delbar,\a\bbar}=u_{\g\delbar \a\bbar}+u_{s\delbar}R_{\a\bbar\g\bar{s}}
.\tag 2.4
$$
Using the fact that
$$
u_{\g\delbar, \a\bar{\a}}=u_{\g\delbar, \bar{\a}\a}+R_{\g\bar{p}}u_{p\delbar}
-R_{p\delbar}u_{\g\bar{p}} \tag 2.5
$$
we have
$$
\split
\Delta u_{\g\delbar}& =\frac{1}{2}(u_{\g\delbar, \a\bar{\a}}+u_{\g\delbar,
\bar{\a}\a})\\
& = u_{\g\delbar \a\bar{\a}}+\frac{1}{2}\left(R_{\g\bar{p}}u_{p\delbar}+
R_{p\delbar}u_{\g\bar{p}}\right)
\endsplit.\tag 2.6
$$
Combining with (2.3), we conclude that $u_\abb$ satisfies (1.2).
\enddemo
In the following, $\tn$ and $\tD$ denote the covariant derivative and the Laplacian with respect with the initial metric.

\proclaim{Proposition 2.2} Let $u_0$ be a smooth function such that
$|u_0(x)|\le  \exp(a(r_0(x)+1))$ for all $x$ for some positive constant $a>0$.
 Let $u(x,t)$ be the solution of (2.1) obtained in Proposition 2.1. We have the following:

\roster
\item"{(i)}" For any $b>0$
$$
\int_0^T \int_M\exp(-br_0^2(x))\lf(|\nabla u|^2(x,t)+t||u_\abb||^2(x,t)\ri)dV_tdt<\infty,\tag2.7
$$
where $||u_\abb||^2=g^{\a\bar\delta}g^{\g\bb}u_\abb u_{\g\bar\delta}.$
\item"{(ii)}" If in addition,
$\int_{B_0(o,r)}|\tn u_0|^2 dV_0\le \exp(a'(1+r))$ for some $a'>0$, where $B_0(o,r)$ is the geodesic ball with center at $o$ and   radius $r$ with respect to the initial metric $g(0)$, then
$$
\int_0^T \int_M\exp(-br_0^2(x)) ||u_\abb||^2(x,t)dV_tdt<\infty,\tag2.8
$$
\item"{(iii)}" If in addition $|\tn u_0|^2\le C_1$ on $M$ then
$$
|\nabla u|^2\le C_1\tag2.9
$$
 and
$$
||u_\abb||^2(x,t)\le \frac{C_2}t\tag2.10
$$
for some constant $C_2$ on $M\times [0,T]$.
\endroster

\endproclaim
\demo{Proof} By Proposition 2.1, there exists positive constant and $c$ such that
$$
|u (x,t)|\le  \exp(c(r_0(x)+1))\tag2.11
$$
on  $M\times [0,T]$. Since
$$
\heat u^2=-|\nabla u|^2,
$$
we can proceed as in the proof of Lemma 1.1 to conclude that for any $b>0$
$$
\int_0^T \int_M\exp(-br_0^2(x))|\nabla u|^2dV_tdt<\infty.
\tag2.12
$$
Direct computations show (see \cite{N-T, Lemma 1.1} for example)
$$
\heat |\nabla u|^2=-||u_{\a\b}||^2-||u_\abb||^2.\tag2.13
$$
Combining with (2.12), one can also proceed as in the proof of Lemma 1.1 and conclude that (2.7) is true. In case $|\tn u_0|^2$ satisfies the condition in (ii),   then one can prove (2.8) similarly.

By (2.13), it is easy to see that $\heat (\sqrt{|\nabla u|^2+1})\le 0$. Suppose $|\tn u|^2\le C_1$ on $M$, then by (2.7) we can apply Theorem 1.1 to conclude that (2.9) is true.

Since $u_\abb$ satisfies (1.2), as in the proof of (1.11) we have
$$
\heat (1+t\Phi)^\frac12  \le C_3\Phi
$$
on $M\times[0,T]$ for some constant $C_3>0$, where $\Phi=||u_\abb||^2$.  Hence on $M\times[0,T]$,
$$
\heat \lf(C_3|\nabla u|^2+ (1+t\Phi)^\frac12\ri)\le0
$$
where we have used (2.13).
By (2.8) and (2.9), we can apply the maximum principle in \cite{N-T} and conclude that
$$
\sup_{M\times [0,T]} \lf(C_3|\nabla u|^2+ (1+t\Phi)^\frac12\ri)\le  C_3C_1+1
$$
where we have used the fact that $|\tn u|^2\le C_1$. From this (2.10) follows.
\enddemo

Next, we shall study the properties of $u(x,t)$ in case the initial value $u_0$ is plurisubharmonic.
\proclaim{Theorem 2.1} Let $u_0(x)$ be a smooth function on $M$ such that  (a) $u_0$ is plurisubharmonic; and (b) there exists $a>0$ such that  $|u_0(x)|\le \exp(a(1+r_0(x))$   and   $\tD u_0\le \exp(a(1+r_0(x))$. Let $u$ be the solution of (2.1) obtained in Proposition 2.1. We have the following:
\roster
\item"{(i)}" $u(x,t)$ is plurisubharmonic for $t>0$.
\item"{(ii)}" If $u_0$ is not harmonic, then $w=u_t>0$ for $t>0$, and
we have the following differential inequality:
$$
w_t-\frac{|\nabla w|^2}w+\frac wt\ge0\tag2.14
$$
for $t>0$.
\endroster
If in addition, $\sup_M|\tn u_0|^2\le C_1<\infty$ for some constant $C_1$, then $u_\abb$ satisfies (2.10) for some constant $C_2$.
\endproclaim
\demo{Proof} Let $f=\tD u_0\ge0$. By assumptions, $|u_0(x)|\le \exp(a(1+r_0(x))$ and $f(x)\le  \exp(a(1+r_0(x))$. It is easy to see that
$$
\int_{B_0(o,r)}|\tn u_0|^2dV_0\le \exp(a'(1+r_0(x))
$$
for some $a'>0$. Hence $u_\abb$ satisfies (2.8) by Proposition 2.2. Since $u_0$ is plurisubharmonic, we also have $||u_\abb||^2(x,0)\le \exp(a''(1+r_0(x))$ for some $a''>0$. By (i), Proposition 1.1 and Lemma 2.1, we conclude that $u$ is plurisubharmonic for $t>0$.

Since $u_\abb$ satisfies (1.2) by Lemma 2.1 and $w=u_t=\D u$, taking trace of (1.2), we have
$$
\heat w=R_\abb u_{\beta\bar \a}\ge0.
$$
If $w(x,t)=0$ for some $x$ and $t>0$, then by the strong maximum principle   (see
\cite{Cw3, Proposition 3.6}), we have $\tD u_0=0$ on $M$. Hence if $u_0$ is not harmonic, then $w>0$ for $t>0$.

Since $u_\abb$ satisfies the conditions in Theorem 1.2,
if we let $h_{\abb}$, in Theorem 1.2
 to be $u_{\abb}$, then in normal coordinates
$$
div(h)_{\a}=\nabla_{\g}u_{\a\bar{\g}}= \nabla_\a (u_t)\ \  \text{ and}\ \
div(h)_{\delbar}=\nabla_{\bar{\a}}u_{\a\delbar}=\nabla_{\delbar}(u_t),
$$
and
$$
Z=\D (u_t) +R_{\a\bbar}u_{\bar{\a}\beta}+\nabla_{\bar{\a}}(u_t)V_{\a}
+\nabla_{\a}(u_t)V_{\bar{\a}}+u_{\a\bbar}V_{\bar{\a}}V_{\beta}+\frac{u_t}{t}
\ge 0
$$
 for any $(1,0)$ vector field $V$. Combining this with (2.15), we have
$$
w_t+\nabla_{\bar{\a}}w V_{\a}
+\nabla_{\a}w V_{\bar{\a}}+u_{\a\bbar}V_{\bar{\a}}V_{\beta}+\frac{w}{t}
\ge 0.
$$
Choosing $V_\a=-\frac{\nabla_{\a} w}{w}$ we conclude that (2.14) is true.

The last assertion follows from Proposition 2.2 immediately.
\enddemo

\proclaim{Remark 2.1} If $u_0(x)$ is a solution to the Poincar\'e-Lelong
equation $\i\ddbar u_0 =\Ric(x,0)$, by Theorem 1.3 of [N-T] we know that
we have a solution $u(x,t)$ to (2.1) in this case with
$\i\ddbar u(x,t)=\Ric(x,t)$. Then (2.14) in Theorem 2.1 is nothing but
the differential LYH inequality
 of Cao on the scalar curvature since $w(x,t)=R(x,t)$.
\endproclaim

Next we shall prove a Li-Yau type differential inequality for the
positive plurisubharmonic solution of (2.1). The result will not be needed in the next section.

\proclaim{Theorem 2.2}
Let $u(x,t)$ be a positive solution to (2.1) such that $u(x,t)$  is
plurisubharmonic for all $t$. Then we have the following differential
inequality:
$$
\frac{u_t}{u}-\frac{|\nabla u|^2}{u^2}+\frac{m}{t}\ge 0 \tag 2.15
$$
\endproclaim
\demo{Proof} As in Li-Yau [L-Y], we let $v=\log u$. Then
 $$\D v -v_t =-|\nabla v|^2. \tag 2.16$$
Let
$$
G(x,t)=t\lf(|\nabla v|^2-\eta v_t\ri),
$$
where $\eta >1$ is a constant.
Direct calculation shows that in   normal coordinates at a point:
$$
\D |\nabla v|^2 =R_{\bar{\a}\beta}v_{\a}v_{\bb}+v_{\a\g}v_{\bar{\a}\bar{\g}}
+v_{\a\bar{\g}}v_{\ba\g} +(\D v)_\a v_{\ba}+v_{\a}(\D v)_{\ba}, \tag 2.17
$$
$$
\frac{\p}{\p t}|\nabla v|^2 =R_{\bar{\a}\beta}v_{\a}v_{\bb}+(v_t)_\a v_{\ba}
+v_{\a}(v_t)_{\ba} \tag 2.18
$$
and
$$
v_{tt}-\D (v_t)=R_{\ba\b}v_{\abb}+R_{\ba\b}v_{\a}v_{\bb}+(v_t)_\a v_{\ba}
+v_{\a} (v_t)_{\ba} \tag 2.19
$$
Combining (2.16)--(2.19) we have that
$$
\split
\lf( \D -\frac{\p}{\p t}\ri)\lf(|\nabla v|^2-\eta v_t\ri) &=
 v_{\a\g}v_{\bar{\a}\bar{\g}}
+v_{\a\bar{\g}}v_{\ba\g}-\lf(|\nabla v|^2-\eta v_t\ri)_{\a}v_{\ba}
-v_{\a}\lf(|\nabla v|^2-\eta v_t\ri)_{\ba}\\
&\ \ +\eta R_{\ba\b}\lf(v_{\abb}+v_{\a}v_{\bb}\ri).
\endsplit
$$
Using the fact that
$$
R_{\ba\b}\lf(v_{\abb}+v_{\a}v_{\bb}\ri)=\frac{1}{u}R_{\ba\b}u_{\abb}\ge 0
$$
we then have
$$
\split
\lf( \D -\frac{\p}{\p t}\ri)G &\ge tv_{\a\bar{\g}}v_{\ba\g}-
2<\nabla G, \nabla v> -\frac{G}{t}\\
& \ge \frac{t}{m}\lf(\D v\ri)^2 -2<\nabla G, \nabla v> -\frac{G}{t}\\
& = \frac{t}{m}\lf(|\nabla v|^2-v_t\ri)^2-2
<\nabla G, \nabla v> -\frac{G}{t}.
\endsplit
\tag 2.20
$$
Once we have (2.20), we can use the cut-off function argument
as in [L-Y] to carry the
interior estimates. For the sake of the completeness we include the argument
here. Let $\psi(s)$ be a cut-off function such that $0\le \psi\le 1$,
$\psi(s)\equiv 1$ for $s \in [0,1]$ and $\psi(s)\equiv 0$ for $s\ge 2$.
We also require that
$$
\psi'\le 0,\ \ \psi'' \ge -C_1 \ \ \text{and }\ \ \frac{|\psi'|^2}{\psi}\le C_1
\tag 2.21
$$
for some positive constant $C_1$. Now we let $\phi(x) =\psi(r_t(x)/R)$.
Let $\Phi=\phi G$. Suppose $\Phi$ attains a positive   maximum at $(x_0,t_0)$. Then  we have at $(x_0,t_0)$.
$$
0\ge \lf(\D -\frac{\p}{\p t}\ri) \Phi\ \  \text{and }\ \  \nabla \Phi =0.
$$
Note that $\phi$ may not be smooth at $x_0$ in the space variable, but we can always use the trick of Calabi as in \cite{L-Y}. $\phi$ may not be  smooth in the $t$ variable at $t_0$, but we can use difference quotient so that the final result of the following computations is correct. Hence the above differential inequality together with (2.20) implies    that at $(x_0,t_0)$
$$
\split
0\ge t_0 \phi \lf(\D -\frac{\p}{\p t}\ri) \Phi & \ge \frac{1}{m}
\lf(t_0\phi |\nabla v|^2 -t_0\phi v_t\ri)^2 -2Gt_0\frac{|\nabla \phi|^2}{\phi}+
t\lf(\D \phi-\frac{\p}{\p t}\phi\ri)G \\
& \ \ -G\phi^2 +2<\nabla \phi, \nabla v>G \phi t_0\\
&\ge \frac{1}{m} \lf(t_0\phi |\nabla v|^2 -t_0\phi v_t\ri)^2-
G\lf[2t_0\frac{|\nabla \phi|^2}{\phi}-t_0\lf(\D \phi-\frac{\p}{\p t} \phi\ri)+1\ri]\\
&\ \ -2\frac{|\nabla \phi|}{\phi^{1/2}}G\lf(|\nabla v|\phi^{1/2}t_0^{1/2}\ri)
t_0^{1/2}.
\endsplit \tag 2.22
$$
Using (2.21) we have that
$$
\frac{|\nabla \phi|^2}{\phi}\le \frac{C_2}{R^2} \ \ \text{ and }\ \
-\D \phi \ge \frac{C_2}{R^2}. \tag 2.23
$$
Also Theorem 17.2 of [H4] implies that
$$
|\frac{\p}{\p t}\phi |\le \frac{C_2}{R}. \tag 2.24
$$
Here $C_2$ is a constant dependent of $C_1$, $m$ and the upper bound of
$|Rm|(x,t)$.
Combining (2.22)--(2.24) we have, at the maximum of $\Phi$ over
$M\times [0, T]$, that
$$
0\ge \frac{1}{m}(y-z)^2-\frac{C_2}{R}(y-\eta z)y^{1/2}t_0^{1/2}-(y-\eta z)\lf(
\frac{C_2t_0}{R^2}+\frac{C_2t_0}{R}+1\ri). \tag 2.25
$$
Here $y=t_0\phi|\nabla v|^2(x_0,t_0)$, $z=t_0\phi v_t(x_0,t_0)$. Using the trick of [L-Y], we
write
$$
(y-z)^2 =\frac{1}{\eta^2}(y-\eta z)^2 +2\frac{\eta-1}{\eta}(y-\eta z)y+
\lf(\frac{\eta-1}{\eta}\ri)^2y^2.
$$
Using the $ax^2+bx \ge -\frac{b^2}{4a}$, for $R>>1$ we have that
$$
0\ge \frac{1}{m\eta^2}(y-\eta z)^2 -(y-\eta z)\lf(\frac{C_3 t_0}{R}+1\ri).
$$
for some constant $C_3$ independent of $R$. From which we have that
$$
\sup_{B_0(o,R)\times[0,T]}t\lf(|\nabla v|^2-\eta v_t\ri) \le m\eta^2\lf( 1+
\frac{C(m,\eta,T, \sup_{M\times [0,T]}{|Rm|})}{R}\ri).
$$
Here we have used the fact that $g(t)$ is nonincreasing so that $B_t(o,R)\supset B_0(o,R)$. Letting $R\to \infty$ and then $\eta\to 1$ we have (2.15).
\enddemo


%% file: HP-s3.tex
\input amstex
\documentstyle{amsppt}
\magnification=1200
\hsize=13.8cm
\def\i{\sqrt{-1}}
\def\Ric{\text{Ric}}
\def\lf{\left}
\def\ri{\right}
\def\bbar{{\bar \beta}}
\def\a{\alpha}
\def\g{\gamma}
\def\p{\partial}
\def\delbar{{\bar\delta}}
\def\ddbar{\partial\bar\partial}

\def\C{\Bbb C}

\def\Rm{{\text{\rm Rm}}}

\def \b {\beta}
\def \e{\epsilon}

\def\ba{{\bar\alpha}}
\def\bb{{\bar\beta}}
\def\abb{{\alpha\bar\beta}}

\def\tD{\widetilde \Delta}
\def\tn{\widetilde \nabla}
\def \D {\Delta}

\def\heat{\lf(\frac{\p}{\p t}-\Delta\ri)}
\def\h{\frak h}
\def\aint{\frac{\ \ }{\ \ }{\hskip -0.4cm}\int}
\subheading{\S3 Liouville properties of plurisubharmonic functions}

In this section, we shall discuss Liouville properties of plurisubharmonic functions using the LYH  type inequality in \S1 and the  results of \S2. In this section, we always assume that $(M,g_\abb(x))$ is a complete noncompact
K\"ahler manifold with bounded nonnegative holomorphic bisectional curvature. We also assume that for all $x\in M$ and $r>0$,  $k(x,r)\le \epsilon(r)$ for some nonincreasing function $\epsilon(r)$ with $\lim_{r\to\infty}\e(r)=0$, where
$$
k(x,r)=\aint_{B_0(x,r)}\Cal R_0dV_0\tag3.1
$$
and $\Cal R_0$ is the scalar curvature of $M,g_\abb(x)$. By \cite{N-T}, we know that  (1.1) has a solution $g_\abb(x,t)$ on $M\times [0,\infty)$ such that for any $0<T<\infty$, $g_\abb$ satisfies (i)-(iv) in \S1.

 Define
$$
F(x,t)=\log\left(\frac{\det(g_{\a\bbar}(x,t))}
{\det(g_{\a\bbar}(x,0))}\right).\tag3.2$$

To illustrate the idea of the proof to a more general result,
let us begin with the following particular case. In this case, what we need is to assume that (1.1) has long time solution $g_\abb$ so that for any $T<\infty$, conditions (i)-(iv) in \S1 are satisfied by $g_\abb$ on $M\times[0,T]$.

\proclaim{Theorem 3.1} With the above assumptions,  suppose $u_0$ is a plurisubharmonic function such that (i) $u$ is bounded; and (ii) $\tD u_0(x)\le \exp(a(1+r_0(x))$ for some constant $a>0$. Then $u_0$ must be constant.
\endproclaim
\demo{Proof} Let $\tD u_0=f$, then $f\ge0$. Since $u_0$ is bounded,
by \cite{N-S-T1, Corollary 2.1} we have
$$
\int_0^\infty s\lf(\aint_{B_0(x,s)}fdV_0\ri)ds\le C_1
$$
for some constant $C_1$ independent of $x$. By \cite{N-S-T1, Corollary 1.2}, we know that
$$\sup_M|\tn u_0 |\le C_2.
$$ By Proposition 2.1, there is a unique solution $u(x,t)$ with initial data $u_0$. Moreover,
 by Proposition 2.1 and the maximum principle in \cite{N-T, Theorem 1.2}, we conclude that $u$ is uniformly bounded.

Since $\tD u_0(x) \le \exp(a(1+r_0(x))$,  by Theorem 2.1(i) we conclude that $u(x,t)$ is plurisubharmonic for all $t>0$. Moreover, suppose $u_0$ is not harmonic, then by Theorem 2.1(ii) $w=u_t>0$ for $t>0$ and $tw$ is nondecreasing in $t$. Hence
$$
\split
u(x,t)-u_0(x)&=\int_0^tw(x,s)ds\\
&\ge w(x,1)\int_1^t\frac 1sds\\
&=w(x,1)\log t.
\endsplit
$$
Since $w(x,1)>0$, let $t\to\infty$, the above inequality contradicts the fact that $u$ is uniformly bounded. Hence $u_0$ must be harmonic and is constant by \cite{Y}.
\enddemo

Next we shall generalize  Theorem 3.1 by relaxing the condition that $u_0$ is bounded. In the following, we  always assume that $u_0$ is a plurisubharmonic function on $M$ such that there exists a constant $a>0$ such that
$$
\cases |\tn u_0(x)|&\le a,\\
\tD u_0(x)&\le \exp(a(1+r_0(x))
\endcases\tag3.3
$$
for all $x\in M$. Note that in the proof of Theorem 3.1, we know that if $u_0$ is bounded, then $u_0$ will satisfies the first inequality of (3.3).

Because of (3.3), let $u$ be the solution of  (2.1) with initial data $u_0$ constructed in Proposition 2.1. By Proposition 2.2, $u$ is plurisubharmonic for all $t\ge0$.  Let $v(x,t)=u(x,t)-u_0(x)$. Also, let $\frak m(t)=\inf_{x\in M}F(x,t)$. Then $\frak m(t)\le 0$, nonincreasing,  and is finite for fixed $t$ by properties (ii) and (iv) in \S1 and the fact that
$F(x,t)=-\int_0^t\Cal R(x,s)ds$, where $\Cal R(x,s)$ is the scalar curvature at time $s$.

\proclaim{Lemma 3.1} With the above assumptions and notations, we have
$$
\tD v-e^Fv_t\ge -\tD u_0.\tag3.4
$$
\endproclaim
\demo{Proof} As in \cite{Shi, p.156}, using the fact that $g_\abb$ is nonincreasing, we have
$$
\tD u\ge e^F\D u=e^Fu_t=e^F v_t.
$$
Hence
$$
\tD v=\tD u-\tD u_0\ge e^Fv_t-\tD u_0.
$$
The result follows.
\enddemo

\proclaim{Lemma 3.2} With the same assumptions and notations as in Lemma 3.1, there is a constant $C$ such that for all $(x,t)\in M\times[0,\infty)$, we have
$$
0\le v(x,t)\le Ct^\frac12\lf(-\frak m(2t)+1\ri)
\tag3.5
$$
\endproclaim
\demo{Proof} First note that $v(x,0)=0$ and $v_t=u_t=\tD u\ge0$.
 Hence $v\ge0$. We need a more refined estimate of (2.10). More precisely, the Bochner formula on $\|u_{\abb}\|^2 $ says that
$$
\lf(\D -\frac{\p}{\p t}\ri)\|u_{\abb}\|^2 \ge \|u_{\abb\g}\|^2+
\|u_{\abb \bar{\g}}\|^2-{\Cal R}(x,t)\|u_{\abb}\|^2.
$$
Using the LYH type inequality of H.-D. Cao  as in [N-T] we have that
$$
t{\Cal R}(x,t)\le -2\frak m(2t).
$$
Combining them we have that
$$
\lf(\D -\frac{\p}{\p t}\ri)\lf(1+t\Phi \ri)^{\frac{1}{2}}
\ge -(-\frak m(2t)+1)\Phi.
$$
Here $\Phi =\|u_{\abb}\|^2$. Now we can proceed as in the proof of
Proposition 2.2 (iii) to conclude that
$$
\sup_{M\times [0,T]}\lf((-\frak m(2T)+1)|\nabla u|^2 +(1+t\Phi)^{\frac{1}{2}}
\ri)\le a^2 (-\frak m(2T)+1) +1,
$$
which then implies
$$
||u_\abb||(x,t)\le C_1t^{-\frac12}\lf(-  \frak m(2t)+1\ri)
$$
for some constant $C_1$ depending only on $m$ and $\sup_M|\tn u_0|$. Hence
$$
\split
v(x,T)&=\int_0^T v_t(x,t)dt\\
&=\int_0^T u_t(x,t)dt\\
&=\int_0^T\D u(x,t)dt\\
&\le C_1\int_0^Tt^{-\frac12}\lf( -\frak m(2t)+1\ri)dt\\
&\le C_2T^\frac12\lf(-\frak m(2T)+1\ri)
\endsplit
$$
for some constant $C_2$ independent of $x$ and $t$. The proof of the lemma is completed.
\enddemo

Using the method of proof of Theorem 2.1 in \cite{N-T}, we have:

\proclaim{Theorem 3.2} Let $u_0$ be a plurisubharmonic function on $M$ satisfying (3.3). Suppose
$$
\limsup_{R\to\infty}\frac{\sup_{\p B_0(o,R)}u_0}{\log t}\le0\tag3.6
$$
where $R^2=t^{\frac 32}e^{-\frak m(2t)} (-\frak m(2t)+1)$, then $u_0$ must be constant.
\endproclaim

\noindent{\bf Remarks}
\roster
\item"{(a)}" It is easy to see that $R\to\infty$ if and only if $t\to\infty$.
\item"{(b)}" By \cite{Sh2-3} and \cite{N-T, Remark 2.2} if $\e(r)$ at the beginning of this section satisfies $\e(r)\le Cr^{-2}$ or more generally if $\int_0^rs\e(s)ds\le C\log (r+2)$, then $-\frak m(t)\le C'\log (t+1)$, and the assumption (3.6) can be replaced by
$$
\limsup_{x\to\infty}\frac{u_0(x)}{\log r_0(x)}\le0.
$$
\item"{(c)}" Similarly, if $\e(r)\le r^{-\theta}$ for some $\theta>0$, then the assumption  (3.6) can be replaced by
$$
\limsup_{x\to\infty}\frac{u_0(x)}{\log\log r_0(x)}\le0.
$$
If $\int_0^rs\e(s)ds\le Cr^2/\log (2+r)$, then the assumption (3.6) can be replaced by
$$
\limsup_{x\to\infty}\frac{u_0(x)}{\log\log\log r_0(x)}\le0,
$$
and so on.

\item"{(d)}" It is easy to see that Theorem 3.1 is a particular case of Theorem 3.2
\endroster

\demo{Proof of Theorem 3.2} Let $u$ and $v$ as in Lemmas 3.1 and 3.2. Let $(x_0,T)\in M\times (0,\infty)$. For any $R>0$, let $G_R$ be the positive Green's function with zero boundary value on $B_0(x_0,R)$ with respect to the initial metric. By (3.4)
$$
\split
\int_0^T&\int_{B_0(x_0,R)}G_R(x_0,y)\tD v(y,t)dV_0dt\\
&\ge -T\int_{B_0(x_0,R)}G_R(x_0,y)\tD u_0(y)dV_0+\int_0^T\int_{B_0(x_0,R)}
G_R(x_0,y) e^{F(y,t)}v_t(y,t)dV_0dt\\
&\ge -T\int_{B_0(x_0,R)}G_R(x_0,y)\tD u_0(y)dV_0+e^{\frak m(T)}\int_{B_0(x_0,R)}G_R(x_0,y)v(y, T)dV_0\\
&\ge C_1(m)\lf[ -T\lf(-u_0(x_0)+\sup_{B_0(x_0,R)}u_0\ri)+e^{\frak m(T)}R^2\aint_{B_0(x_0,\frac R5)}v(y,T)dV_0\ri]
\endsplit\tag3.7
$$
for some positive constant $C_1$ depending only on $m$, where we have used Theorem 2.1 in \cite{N-S-T1} and Lemma 2.2 in \cite{N-T} and  the fact that
$v\ge0$, $F_t\le 0$. On the other hand,   by Green's formula and Lemma 3.2, we have that, for any $0<t<T$,
$$
\split
\int_{B_0(x_0,R)}G_R(x_0,y)\tD v(y)dV_0&=-v(x_0,t)-\int_{\p B_0(x_0,R)}v \frac{\p G_R}{\p \nu}\\
&\le C_2t^\frac12\lf(-\frak m(2t)+1\ri)
\endsplit \tag3.8
$$
for some constant $C_2$ independent of $(x,t)$. By (3.4) and the fact that $v_t\ge0$, we have $\tD v\ge -\tD u_0$. Since $v\ge0$, by the generalized mean value inequality \cite{N-T, Lemma 2.1}, (3.7) and (3.8), we have
$$
\split
v(x_0,T)&\le C_3\aint_{B_0(x_0,\frac R5)}v(y,T) dV_0+\int_{B_0(x_0,\frac R5)}G_{\frac R5}(x_0,y)\tD u_0(y)dV_0\\
&\le C_4\bigg[R^{-2}Te^{-\frak m(2T)}\lf( -u_0(x_0)+\sup_{B_0(o,2R)}u_0+T^\frac12 (-\frak m(2T)+1)\ri)\\
&-u_0(x_0)+\sup_{B_0(o, 2R)}u_0\bigg]
\endsplit
$$
if $R$ is large, for some constants $C_3$ and $C_4$ independent of $(x_0,T)$ and $R$. Let $R$ be such that $(2R)^2=T^{\frac 32}(1+T)e^{-\frak m(2T)} (-\frak m(2T)+1)$, then by (3.6), we can conclude that
$$
\limsup_{t\to\infty}\frac{v(x_0,t)}{\log t}=0.\tag3.9
$$
We claim that $u_0$ is harmonic. Suppose not, then as in the proof of Theorem 3.1, we have $u(x_0,t)\ge C\log t$ for some   constant $C>0$ for all $t\ge 1$. This is impossible.

By the definition of $R$ in (3.6), it is easy to see that $\log R\ge \log t$ when $t$ is large. Hence (3.6)  implies that
$$
\limsup_{R\to\infty}\frac{\sup_{\p B_0(o,R)}u_0}{\log R}=0.
$$
Since $u_0$ is harmonic, it must be constant by \cite{C-Y}.
\enddemo

Since one can solve the Poincar\'e-Lelong equation for a (1,1) form on a complete noncompact manifold with nonnegative holomorphic bisectional curvature under rather weak assumptions on the (1,1) form  (see \cite{N-S-T1}), one can apply Theorem 3.2 (or Theorem 3.1) to obtain results on  the flatness of the holomorphic line bundles. As an example, we have the following:

\proclaim{Corollary 3.1}
Let $(M, g_{\abb}(x))$ be a complete nocompact K\"ahler manifold with bounded
nonnegative holomorphic bisectional curvature satisfying the conditions in Theorem 3.1. Let $(L, \h_0)$ be a holomorphic line bundle
on $M$ with the Hermitian metric $\h_0$. Suppose     $\Omega(\h_0)\ge 0$ and suppose its trace
${\Cal S}_0=g^{\abb}(x)\Omega_{\abb}(\h_0)(x)$ is bounded and
$$
  \int_0^\infty s\aint_{B_0(x,s)}{\Cal S}_0(y)\, dy\, ds\le C\tag3.10
$$
for some  constant $C>0$  for all $x\in M$. Then $(L, \h_0)$ is flat.
\endproclaim
\demo{Proof} Using the fact that $\Cal S_0$ is bounded and (3.10), one can find bounded function $u_0$ such that $\i\ddbar u_0=\Omega(\h_0)$    by \cite{N-S-T1, Theorem 5.1}. Since $\Omega(\h_0)$ is nonnegative, $u_0$ is plurisubharmonic. By Theorem 3.1, $u_0$ is constant and hence $(L, \h_0)$ is flat.
\enddemo
